\documentclass[11pt]{article}

\usepackage{color}
\usepackage{latexsym}
\usepackage{dsfont}
\usepackage{amssymb}
\usepackage{graphicx}
\usepackage{amsmath, amsfonts,amssymb,theorem,euscript,array,enumerate,amsfonts,mathrsfs}
\usepackage{hyperref}

\newtheorem{Theorem}{Theorem}[part]

\newtheorem{Proposition}{Proposition}[part]

\newtheorem{Lemma}{Lemma}[part]
\newtheorem{Corollary}{Corollary}[part]
\newtheorem{Remark}{Remark}[part]

\def \trans{^{\scriptscriptstyle{\intercal}}}

\def \Frac{\displaystyle\frac}

\def \b1{\bf{1}}

\def \N{\mathbb{N}}
\def \R{\mathbb{R}}
\def \S{\mathbb{S}}

\def \E{\mathbb{E}}
\def \F{\mathbb{F}}
\def \G{\mathbb{G}}
\def \P{\mathbb{P}}

\def \T{\mathbb{T}}

\def\esssup_#1{\underset{#1}{\mathrm{ess\,sup\, }}}

\def\argmin_#1{\underset{#1}{\mathrm{argmin\, }}}

\def \Ac{{\cal A}}
\def \Bc{{\cal B}}

\def \Ec{{\cal E}}
\def \Fc{{\cal F}}
\def \Gc{{\cal G}}

\def \Kc{{\cal K}}
\def \Lc{{\cal L}}
\def \Pc{{\cal P}}

\def \Sc{{\cal S}}

\def \Uc{{\cal U}}

\def \Yc{{\cal Y}}

\def \Zc{{\cal Z}}

\def \ni{\noindent}

\def \eps{\varepsilon}

\def \ep{\hbox{ }\hfill$\Box$}

\def\Dt#1{\Frac{\partial #1}{\partial t}}

\def\reff#1{{\rm(\ref{#1})}}

\def\beqs{\begin{eqnarray*}}
\def\enqs{\end{eqnarray*}}
\def\beq{\begin{eqnarray}}
\def\enq{\end{eqnarray}}

\begin{document}

\title{
Discrete time approximation of  fully nonlinear HJB equations via  BSDEs with nonpositive jumps 
}


\author{Idris Kharroubi\footnote{The research of the author benefited from the support of the French ANR research grant LIQUIRISK (ANR-11-JS01-0007).}\\ \footnotesize{CEREMADE, CNRS UMR 7534}, \\\footnotesize{Universit\'e Paris Dauphine}\\ \footnotesize{and CREST,} \\
\footnotesize{ \texttt{kharroubi at  ceremade.dauphine.fr}} 
  \and Nicolas Langren\'e\\ \footnotesize{Laboratoire de Probabilit\'es et Mod\`eles Al\'eatoires,}\\ \footnotesize{Universit\'e Paris Diderot}\\ \footnotesize{ and EDF R\&D}\\\footnotesize{\texttt{langrene at math.univ-paris-diderot.fr}}
   \and Huy\^en Pham\\ \footnotesize{Laboratoire de Probabilit\'es et Mod\`eles Al\'eatoires,}\\ \footnotesize{Universit\'e Paris Diderot}\\ \footnotesize{and CREST-ENSAE}\\\footnotesize{\texttt{pham  at math.univ-paris-diderot.fr}}
}

\maketitle 
%

\begin{abstract}
We propose a new probabilistic numerical scheme for fully nonlinear equation of Hamilton-Jacobi-Bellman (HJB) type associated to stochastic control problem,  which is based  on  the Feynman-Kac representation in \cite{kp12} by means of control randomization and  backward stochastic differential equation with nonpositive jumps.  We study a discrete time approximation for the mini\-mal solution to this class of BSDE when the time step goes to zero, 
which provides both an approximation for the value function and for an optimal control in feedback form. 
We obtained a convergence rate  without any ellipticity condition on the  controlled diffusion coefficient. Explicit implementable scheme based on Monte-Carlo simulations and empirical regressions,   associated error analysis, and numerical experiments are performed in the companion paper \cite{khalanpha13b}. 
\end{abstract}

\vspace{5mm}


\vspace{5mm}


\newpage


\section{Introduction}

Let us consider the  fully nonlinear generalized Hamilton-Jacobi-Bellman (HJB) equation:
\begin{equation}\label{HJBintro} 
\left\{
\begin{array}{rcl}
\Dt{v} + \sup_{a \in A} \big [ b(x,a). D_x v + \frac{1}{2} {\rm tr}(\sigma\sigma\trans(x,a)D_x^2 v) + f(x,a,v,\sigma\trans(x,a)D_xv) \big] &=& 0, \;  \mbox{ on } [0,T) \times \R^d, \\
v(T,x) &=& g, \;  \mbox{ on } \R^d.
\end{array}
\right.
\end{equation}
In  the particular case where $f(x,a)$ does not depend on $v$ and $D_x v$, this partial differential equation (PDE)  is the dynamic programming equation for the stochastic control problem:
\beq \label{constov}
v(t,x) &=& \sup_\alpha \E \Big[ \int_t^T f(X_s^{\alpha},\alpha_s) ds + g(X_T^{\alpha}) \Big| X_t^\alpha = x \Big],
\enq
with controlled diffusion in $\R^d$: 
\beqs
dX_t^\alpha &=& b(X_t^\alpha,\alpha_t) dt + \sigma(X_t^\alpha,\alpha_t) dW_t,
\enqs 
and where $\alpha$ is an adapted control process  valued in a compact space $A$ of $\R^q$. Numerical methods for  parabolic partial differential equations (PDEs) are largely developed in the literature, 
but remain a  big challenge for fully nonlinear PDEs, like the HJB equation \reff{HJBintro}, especially in high dimensional cases. We refer to the recent paper \cite{ftw11} for a review of some deterministic and probabilistic approaches.

In this paper, we propose a new probabilistic numerical scheme for HJB equation, relying on the following Feynman-Kac formula  for HJB equation obtained by randomization of the control process $\alpha$.  
We consider the minimal solution $(Y,Z,U,K)$ to the backward stochastic differential equation (BSDE)  with nonpositive jumps: 
\begin{equation} \label{BSDEnonpos}
\left\{ \begin{array}{ccl}
Y_t &=& g(X_T) +  \int_t^T f(X_s,I_s,Y_s,Z_s) ds + K_T - K_t \\
& & \;\; - \int_t^T Z_s dW_s - \int_t^T \int_A U_s(a) \tilde\mu(ds,da),  \;\;\;\;\;   0 \leq t \leq T, \\
U_t(a) & \leq & 0,  
\end{array}
\right.
\end{equation}
with a forward Markov regime-switching diffusion process $(X,I)$ valued in $\R^d\times A$ given by:
\beqs
X_t &=& X_0 + \int_0^t b(X_s,I_s) ds + \int_0^t \sigma(X_s,I_s) dW_s \\
I_t &=& I_0 +  \int_{(0,t]} \int_A (a-I_{s^-}) \mu(ds,da).
\enqs
Here $W$ is a standard Brownian motion, $\mu(dt,da)$ is a Poisson random measure  on $[0,\infty)\times A$ with finite intensity measure $\lambda(da)$ of full topological support on  $A$, 
and compensated measure $\tilde\mu(dt,da)$ $=$ $\mu(dt,da)-\lambda(da)dt$. Assumptions on the coefficients $b,\sigma,f,g$ will be detailed in the next section, but  we emphasize the important point 
that no degeneracy condition on the controlled diffusion coefficient $\sigma$ is imposed.  
It is proved in \cite{kp12} that the minimal solution  to this class of BSDE is related  to the HJB equation  \reff{HJBintro}  through the relation $Y_t$ $=$ $v(t,X_t)$.

The purpose of this paper is to provide and analyze a discrete-time approximation scheme for the minimal solution to \reff{BSDEnonpos}, and thus an approximation scheme for the HJB equation. In the non-constrained jump case, approximations schemes for BSDE have been studied in  the papers \cite{gobetal06}, \cite{BE08}, which extended works in   \cite{boutou04}, \cite{zha04} for BSDEs in a Brownian framework.  The issue is now to deal with the nonpositive jump constraint in \reff{BSDEnonpos}, and we propose a discrete time approximation scheme of the form: 
\begin{equation}\label{schemeBSDE}
\left\{ \begin{array}{rcl}
\bar Y_T^\pi \; = \; \bar \Yc_T^\pi &=& g(\bar X_T^\pi) \\
\bar \Zc_{t_k}^\pi &=& \E \Big[ \bar Y_{t_{k+1}}^\pi \frac{W_{t_{k+1}} - W_{t_k}}{t_{k+1}-t_k} \big| \Fc_{t_k} \Big] \\
\bar\Yc_{t_k}^\pi &=& \E\Big[ \bar Y_{t_{k+1}}^\pi  \big| \Fc_{t_k} \Big] +  (t_{k+1}-t_k) \; f(\bar X_{t_k}^\pi,I_{t_k},\bar\Yc_{t_k}^\pi,\bar \Zc_{t_k}^\pi) \\
\bar Y_{t_k}^\pi &=& \esssup_{a \in A} \E \Big[ \bar\Yc_{t_k}^\pi \big| \Fc_{t_k},  I_{t_k} = a \Big], \;\;\; k = 0, \ldots, n-1, 
\end{array}
\right.
\end{equation}
where $\pi$ $=$ $\{t_0=0 < \ldots <t_k < \ldots <t_n = T\}$ is a partition of the time interval $[0,T]$,  with modulus $|\pi|$, and $\bar X^\pi$ is the Euler scheme of $X$ (notice that $I$ is perfectly simulatable once we know how to simulate the distribution 
$\lambda(da)/\int_A \lambda(da)$ of the jump marks). The interpretation of this scheme is the following. The 
first three lines in \reff{schemeBSDE} correspond to the standard scheme $(\bar \Yc^\pi,\bar \Zc^\pi)$ for a discretization of a BSDE with jumps (see \cite{BE08}), where we omit here the computation of the jump component. The last 
line in \reff{schemeBSDE} for computing the approximation $\bar Y^\pi$ of the minimal solution $Y$  corresponds precisely to the minimality condition for the nonpositive jump constraint and should be understood as follows. By the Markov property of the forward process $(X,I)$, the solution $(\Yc,\Zc,\Uc)$  to the BSDE with jumps (without constraint)  is in the form $\Yc_t$ $=$ $\vartheta(t,X_t,I_t)$ for some deterministic function $\vartheta$.  Assuming that $\vartheta$ is a continuous function, 
the jump component of the BSDE, which is induced by a jump of the forward component $I$,   is equal to  $\Uc_t(a)$ $=$ $\vartheta(t,X_t,a)-\vartheta(t,X_t,I_{t^-})$.  Therefore, the nonpositive jump constraint  means that: 
$\vartheta(t,X_t,I_{t^-})$ $\geq$ $\esssup_{a\in A} \vartheta(t,X_t,a)$. The  minimality condition is thus written as: 
\beqs
Y_t \; = \;  v(t,X_t) \; = \;  \esssup_{a \in A} \vartheta(t,X_t,a) &=& \esssup_{a \in A} \E [ \Yc_t | X_t, I_{t} = a ],
\enqs
whose discrete time version is the last line in scheme \reff{schemeBSDE}.  

In this work, we mainly consider  the case where $f(x,a,y)$ does not depend on $z$, and our aim is
to analyze the discrete time approximation error on $Y$, where we split the error between the positive and negative parts:
\beqs
{\rm Err}_+^\pi(Y) \;:=  \;  \Big(\max_{k \leq n-1}  \E \Big[ \big(Y_{t_k} - \bar Y_{t_{k}}^\pi \big)_+^2  \Big]  \Big)^{1\over 2}, & &
{\rm Err}_-^\pi(Y) \; := \;  \Big(\max_{k \leq n-1}  \E \Big[ \big(Y_{t_k} - \bar Y_{t_{k}}^\pi \big)_-^2  \Big]  \Big)^{1\over 2}. 
\enqs
We do not study directly the error on $Z$, and 
instead focus on the approximation of an op\-timal control for the HJB equation, which is more relevant in practice. It appears that the maximization step in the scheme \reff{schemeBSDE} provides a  control in feedback form  $\{\hat a(t_k,\bar X_{t_k}^\pi)$, $k$ $\leq$ $n-1\}$, which approximates the optimal control  with an estimated error bound.  
The analysis of the error on $Y$ proceeds as follows. We first introduce the solution $(Y^\pi,\Yc^\pi,\Zc^\pi,\Uc^\pi)$ of a discretely jump-constrained BSDE. This corresponds formally to BSDEs for which the nonpositive jump constraint operates only a finite set of times, and should be viewed as the analog of  discretely reflected BSDEs  defined in \cite{balpag03} and \cite{boucha08} in the context of the approximation for reflected BSDEs.  By combining BSDE methods and PDE approach with comparison principles, and further with the shaking coefficients method of Krylov \cite{kr00} and Barles, Jacobsen \cite{barjac07}, we prove the monotone convergence of this discretely jump-constrained BSDE  towards the minimal solution to the BSDE with nonpositive jump constraint,  and obtained a convergence rate without any ellipticity condition on the diffusion coefficient $\sigma$.   
We next focus on the approximation error between the discrete time scheme in \reff{schemeBSDE} and the discretely jump-constrained BSDE. The standard argument for studying rate of convergence of such error consists in getting an estimate of the error at time $t_k$: $\E[ |Y_{t_k}^\pi - \bar Y_{t_k}^\pi|^2]$ in function of the same estimate at time $t_{k+1}$, and then conclude by induction 
together with classical estimates for the forward Euler scheme. 
 However, due to the supremum in the conditional expectation in the scheme \reff{schemeBSDE} for passing from $\bar\Yc^\pi$ to $\bar Y^\pi$, 
 which is a nonlinear operation violating the law of iterated conditional expectations,  such argument does not work anymore. Instead, we consider the auxiliary error control  at time $t_k$:
\beqs
\Ec_k^\pi(\Yc) &:=& \E\Big[\esssup_{a\in A}\E_{t_1,a}\big[ \ldots \esssup_{a\in A} \E_{t_k,a}\big[ |\Yc_{t_k}^\pi - \bar \Yc_{t_k}^\pi|^2 \big] \ldots \big] \Big],
\enqs
where $\E_{t_k,a}[.]$ denotes the conditional expectation $\E[.|\Fc_{t_k},I_{t_k}=a]$, and we are able to express $\Ec_k^\pi(\Yc)$ in function of 
$\Ec_{k+1}^\pi(\Yc)$. We define similarly an error control  $\Ec_k^\pi(X)$ for the forward Euler scheme, and prove that it converges to zero with a rate  
$|\pi|$. 
Proceeding by induction,  we then obtain a rate of convergence 
$|\pi|$ for $\Ec_k^\pi(\Yc)$, and consequently  for  $\E[ |Y_{t_k}^\pi - \bar Y_{t_k}^\pi|^2]$.  This leads finally to a 
rate $|\pi|^{1\over 2}$ for ${\rm Err}_-^\pi(Y)$,   $|\pi|^{1\over 10}$ for ${\rm Err}_+^\pi(Y)$, and so $|\pi|^{1\over 10}$ for the global error 
${\rm Err}^\pi(Y)$ $=$ ${\rm Err}_+^\pi(Y)$ $+$ ${\rm Err}_-^\pi(Y)$. 
Moreover, in the  case where $f(x,a)$ does not depend on $y$ (i.e. the case of  standard HJB equation and stochastic control problem),  we obtain a better rate of order $|\pi|^{1\over 6}$ by relying on a stochastic control representation of the discretely jump-constrained BSDE, and by using a convergence rate result in \cite{kr99} for the approximation of controlled diffusion by means of  piece-wise constant policies. 
 Anyway,  our  result improves the convergence rate of the mixed Monte-Carlo finite difference scheme proposed in \cite{ftw11}, where the authors obtained a rate $|\pi|^{1\over 4}$ on one side and  $|\pi|^{1\over 10}$ on the other side under  a nondegeneracy condition.  

We conclude this introduction by pointing out that the above discrete time scheme is not yet directly implemented in practice, and requires the estimation and computation of the conditional expectations together with the supremum. Actually, simulation-regression methods  on  basis functions defined on $\R^d\times A$ appear to be very efficient, and provide  approximate optimal controls in feedback forms via the maximization operation in the 
last step of the scheme \reff{schemeBSDE}. We postpone this analysis and illustrations with several numerical tests arising in superreplication of options under uncertain volatility and correlation in a companion paper \cite{khalanpha13b}.  Notice that since it relies  on the simulation of the forward process $(X,I)$, our scheme does not suffer the  curse of dimensionality encountered in finite difference scheme or controlled Markov chains methods (see \cite{kusdup92},  \cite{bonzid03}), and takes advantage of the high-dimensional properties of Monte-Carlo methods.

The remainder of the paper is organized as follows. In Section 2, we state some useful  auxiliary error estimate for the Euler scheme of the regime switching forward process. We introduce in Section 3 discretely jump-constrained BSDE and relate it to a system of integro-partial differential equations. 
Section 4 is devoted to the convergence of discretely jump-constrained BSDE to the minimal solution of BSDE with nonpositive jumps. 
We provide in Section 5 the approximation error for our discrete time scheme,  and as a byproduct an estimate for the approximate optimal control in the case of classical HJB equation associated to stochastic control problem.

\section{The forward regime switching process}

\setcounter{equation}{0} \setcounter{Assumption}{0}
\setcounter{Theorem}{0} \setcounter{Proposition}{0}
\setcounter{Corollary}{0} \setcounter{Lemma}{0}
\setcounter{Definition}{0} \setcounter{Remark}{0}

Let $(\Omega,\Fc,\P)$ be a probability space supporting  $d$-dimensional Brownian motion $W$, and a Poisson random measure $\mu(dt,da)$ with intensity measure $\lambda(da)dt$ on $[0,\infty)\times A$, where $A$ is a compact set of  $\R^q$, endowed with its Borel tribe $\Bc(A)$, and $\lambda$ is a finite measure on $(A,\Bc(A))$ with full topological support. 
 We denote by $\F$ $=$ $(\Fc_t)_{t\geq 0}$ the completion of the natural filtration generated by $(W,\mu)$, and by $\Pc$ the  $\sigma$-algebra of $\F$-predictable subsets of $\Omega\times\R_+$. 

We fix a finite time horizon $T$ $>$ $0$, and consider the solution $(X,I)$ on $[0,T]$ of the regime-switching diffusion model:
\begin{equation}
\label{forwardXI}
\left\{
\begin{array}{ccl}
X_t &=& X_0 + \int_0^t b(X_s,I_s) ds + \int_0^t \sigma(X_s,I_s) dW_s \\
I_t &=& I_0 +  \int_{(0,t]} \int_A (a-I_{s^-}) \mu(ds,da),
\end{array}
\right.
\end{equation}
where $(X_0,I_0)$ $\in$ $\R^d\times A$,  $b$ $:$ $\R^d\times A$ $\rightarrow$ $\R^d$ and $\sigma$ $:$ $\R^d\times A$ $\rightarrow$ $\R^{d\times d}$,  
are measurable functions, satisfying the Lipschitz  condition: 

\vspace{2mm}

\ni \textbf{(H1)}
There exists a constant $L_1$ such that 
\beqs
|b(x,a)-b(x',a')|+|\sigma(x,a)-\sigma(x',a')| & \leq & L_1 \big(|x-x'|+|a-a'|\big), 
\enqs
for all $x,x'\in \R^d$ and  $a,a'\in A$.  The assumption {\bf (H1)} stands in force throughout the paper, and in this section, 
we shall denote by $C_1$ a generic positive constant  which depends only on $L_1$, $T$, $(X_0,I_0)$ 
and  $\lambda(A)$ $<$ $\infty$, and may vary from lines to lines. Under  {\bf (H1)}, we have the existence and uniqueness of a solution to \reff{forwardXI}, and in the sequel, we 
shall denote by  $(X^{t,x,a},I^{t,a})$  the solution to \reff{forwardXI} starting from $(x,a)$ at time $t$. 

\begin{Remark}
{\rm We do not make any ellipticity assumption on $\sigma$. In particular, some lines and columns of $\sigma$ may be equal to zero, and so there is no loss of generality by considering that the dimension $d$ of $X$ and $W$ are equal. 
\ep
}
\end{Remark}

\vspace{2mm}

We first study the discrete-time approximation of the forward process. Denoting by $(T_n,\iota_n)_n$ the jump times and marks associated to $\mu$, we 
observe that $I$ is explicitly written as:
\beqs
I_t &=& I_0 \mathds{1}_{[0,T_1)}(t) + \sum_{n\geq 1} \iota_n \mathds{1}_{[T_n,T_{n+1})}(t), \;\;\; 0 \leq t \leq T,
\enqs
where the jump times $(T_n)_n$ evolve according  to a Poisson distribution of parameter $\lambda$ $:=$ $\int_A \lambda(da)$ $<$ $\infty$, and the 
i.i.d. marks $(\iota_n)_n$ follow a probability distribution  $\bar\lambda(da)$ $:=$ $\lambda(da)/\lambda$.  Assuming that one can simulate the probability distribution $\bar\lambda$, we then see that  the pure jump process $I$ is perfectly simulated.  
Given a partition $\pi$ $=$ $\{t_0=0 < \ldots <t_k < \ldots t_n = T\}$ of  $[0,T]$,  we shall use the natural Euler scheme $\bar X^\pi$ for $X$, defined by: 
\beqs
\bar X_0^\pi &=& X_0 \\
\bar X_{t_{k+1}}^\pi &=& \bar X_{t_k}^\pi + b(\bar X_{t_k}^\pi,I_{t_k}) (t_{k+1}-t_k) + \sigma(\bar X_{t_k}^\pi,I_{t_k}) (W_{t_{k+1}}-W_{t_k}),  
\enqs
for $k$ $=$ $0,\ldots,n-1$. We denote as usual by  $|\pi|$ $=$ $\max_{k\leq n-1}(t_{k+1}-t_k)$ the modulus of $\pi$, and assume that $n|\pi|$ is bounded 
by a constant independent of $n$, which holds for instance when the grid is regular, i.e. $(t_{k+1}-t_k)$ $=$ $|\pi|$ for all $k$ $\leq$ $n-1$. We also define 
the continuous-time version of $\bar X^\pi$ by setting: 
\beqs
\bar X_t^\pi &=& \bar X_{t_k}^\pi + b(\bar X_{t_k}^\pi,I_{t_k}) (t -t_k) + \sigma(\bar X_{t_k}^\pi,I_{t_k}) (W_{t}-W_{t_k}), \;\; t \in [t_k,t_{k+1}], \; k < n. 
\enqs

By standard arguments, see e.g. \cite{klopla92}, one can obtain  under {\bf (H1)}  the $L^2$-error estimate for the above Euler scheme: 
\beqs
\E \Big[ \sup_{t \in [t_k,t_{k+1}]} \big| X_{t} - \bar X_{t_k}^\pi\big|^2 \Big] & \leq & C_1 |\pi|, \;\;\; k < n. 
\enqs
For our purpose, we shall need a stronger result, and introduce the following error control for the Euler scheme: 
\beq \label{errorauxX}
\Ec_k^\pi(X) &:=&  \E\Big[\esssup_{a\in A}\E_{t_1,a}\big[ \ldots \esssup_{a\in A}  
\E_{t_k,a}\big[ \sup_{t \in [t_k,t_{k+1}]} |X_t - \bar X_{t_k}^\pi|^2 \big] \ldots \big] \Big],
\enq
where $\E_{t_k,a}[.]$ denotes the conditional expectation $\E[.|\Fc_{t_k},I_{t_k}=a]$.  We also denote by $\E_{t_k}[.]$ the 
conditional expectation $\E[.|\Fc_{t_k}]$.  Since $I_{t_k}$ is $\Fc_{t_k}$-measurable, and by the law of 
iterated conditional expectations, we notice that
\beqs
\E \Big[ \sup_{t \in [t_k,t_{k+1}]} \big| X_{t} - \bar X_{t_k}^\pi\big|^2 \Big] & \leq & \Ec_k^\pi(X), \;\;\;  k < n. 
\enqs

\begin{Lemma} \label{lemEulerX}
We have 
\beqs
\max_{k < n} \Ec_k^\pi(X) & \leq & C_1 |\pi|. 
\enqs
\end{Lemma}
{\bf Proof.}  From the definition of the Euler scheme, and under the growth linear condition in {\bf (H1)}, we easily see that 
\beq \label{estimXbar}
\E_{t_k}\Big[   \big| \bar X_{t_{k+1}}^\pi\big|^2 \Big] & \leq & C_1 \big(1 +  \big| \bar X_{t_{k}}^\pi\big|^2 \big), \;\; k < n. 
\enq
From the definition of the continuous-time Euler scheme, and by Burkholder-Davis-Gundy inequality, it is also clear that
\beq \label{estimXbarpi}
\E_{t_k}\Big[ \sup_{t\in[t_k,t_{k+1}]}\big| \bar X_t^\pi - \bar X_{t_k}^\pi \big|^2 \Big] & \leq & C_1( 1 +  \big| \bar X_{t_{k}}^\pi\big|^2 \big) |\pi|, \;\; k < n. 
\enq
We also have the standard estimate for the pure jump process $I$ (recall that $A$ is assumed to be compact and $\lambda(A)$ $<$ $\infty$): 
\beq \label{estimI}
\E_{t_k}\Big[ \sup_{t\in[t_k,t_{k+1}]}\big| I_s - I_{t_k} \big|^2 \Big] & \leq & C_1 |\pi|. 
\enq
Let us denote by $\Delta X_t$ $=$ $X_t-\bar X_t^\pi$, and apply It\^o's formula to $|\Delta X_t|^2$ so that for all $t$ $\in$ $[t_k,t_{k+1}]$:
\beqs
|\Delta X_t|^2 &=& |\Delta X_{t_k}|^2 +  \int_{t_k}^t 2 \big(b(X_s,I_s) - b(\bar X_{t_k}^\pi,I_{t_k})\big). \Delta X_s  +  
\big|\sigma(X_s,I_s) - \sigma(\bar X_{t_k}^\pi,I_{t_k})\big|^2 ds \\
& & \; + \; 2 \int_{t_k}^t (\Delta X_s)' \big( \sigma(X_s,I_s) - \sigma(\bar X_{t_k}^\pi,I_{t_k}) \big) dW_s \\
& \leq &   |\Delta X_{t_k}|^2 +  C_1  \int_{t_k}^{t}  |\Delta X_s|^2 +  |\bar X_s^\pi - \bar X_{t_k}^\pi|^2 + |I_s - I_{t_k}|^2 ds  \\
& & \; + \; 2     \int_{t_k}^t (\Delta X_s)' \big( \sigma(X_s,I_s) - \sigma(\bar X_{t_k}^\pi,I_{t_k}) \big) dW_s, 
\enqs
from the Lipschitz condition on $b$, $\sigma$ in {\bf (H1)}. By taking conditional expectation in the above inequality, we then get:
\beqs
\E_{t_k}  \Big[   |\Delta X_t|^2 \Big] & \leq &  |\Delta X_{t_k}|^2  
+  C_1 \int_{t_k}^{t} \E_{t_k} \big[ |\Delta X_s|^2 + |\bar X_s^\pi - \bar X_{t_k}^\pi|^2 + |I_s - I_{t_k}|^2 \big] ds   \\
& \leq & |\Delta X_{t_k}|^2  + C_1( 1 +  \big| \bar X_{t_{k}}^\pi\big|^2 \big)  |\pi|^2 + C_1 \int_{t_k}^t \E_{t_k} \big[ |\Delta X_s|^2 \big]  ds, \;\;\; t \in [t_k,t_{k+1}], 
\enqs
by \reff{estimXbarpi}-\reff{estimI}.  From Gronwall's lemma, we thus deduce that 
\beq \label{estimXk}
\E_{t_k}  \Big[   |\Delta X_{t_{k+1}}|^2 \Big] & \leq & e^{C_1|\pi|} |\Delta X_{t_k}|^2  + C_1 ( 1 +  \big| \bar X_{t_{k}}^\pi\big|^2 \big)  |\pi|^2 , \;\;\; k < n. 
\enq
Since the right hand side of \reff{estimXk} does not depend on $I_{t_k}$, this shows that
\beqs
\esssup_{a \in A} \E_{t_k,a}  \Big[   |\Delta X_{t_{k+1}}|^2 \Big] & \leq & e^{C_1|\pi|}  |\Delta X_{t_k}|^2  + C_1 ( 1 +  \big| \bar X_{t_{k}}^\pi\big|^2 \big)  
|\pi|^2. 
\enqs
By taking conditional expectation w.r.t. $\Fc_{t_{k-1}}$ in the above inequality, using again estimate \reff{estimXk}  together with \reff{estimXbar} at step $k-1$, and iterating this backward procedure until the initial time $t_0$ $=$ $0$, we obtain: 
\beq
 & & \E\Big[\esssup_{a\in A} \E_{t_1,a}\big[ \ldots \esssup_{a\in A}  
\E_{t_k,a}\big[ |\Delta X_{t_{k+1}}|^2  \big] \ldots \big] \Big] \nonumber \\
& \leq &  e^{C_1n|\pi|} |\Delta X_0|^2  + C_1(1+|X_0|^2) |\pi|^2    \frac{e^{C_1n|\pi|} -1}{e^{C_1|\pi|} - 1 } \nonumber \\
& \leq & C_1 |\pi|, \label{interXpi}
\enq
since $\Delta X_0$ $=$ $0$ and $n|\pi|$ is bounded. 

Moreover, the process $X$ satisfies the standard conditional estimate similarly as for the Euler scheme: 
\beqs
\E_{t_k}\Big[   \big| X_{t_{k+1}} \big|^2 \Big] & \leq & C_1 \big(1 +  \big| X_{t_{k}} \big|^2 \big),  \\
\E_{t_k}\Big[ \sup_{t\in[t_k,t_{k+1}]}\big| X_t  -  X_{t_k} \big|^2 \Big] & \leq & C_1( 1 +  \big| X_{t_{k}}\big|^2 \big) |\pi|, \;\; k < n,
\enqs
from which we deduce by backward induction on the conditional expectations: 
\beq \label{interXpi2}
 \E\Big[\esssup_{a\in A} \E_{t_1,a}\big[ \ldots \esssup_{a\in A}  
\E_{t_k,a}\big[ \sup_{t\in[t_k,t_{k+1}]}\big| X_t  -  X_{t_k} \big|^2 \big] \ldots \big] \Big]  & \leq & C_1 |\pi|. 
\enq  
Finally, by writing that $\sup_{t\in[t_k,t_{k+1}]}|X_t-\bar X_{t_k}^\pi|^2$ $\leq$ $2\sup_{t\in[t_k,t_{k+1}]}|X_t- X_{t_k}|^2$ $+$ $2\Delta X_{t_k}$,  taking successive condition expectations w.r.t to $\Fc_{t_\ell}$ and  essential supremum over $I_{t_\ell}$ $=$ $a$,  for $\ell$ going recursively from $k$ to $0$, we get:
\beqs
\E_{t_k} \Big[ \sup_{t\in[t_k,t_{k+1}]}|X_t-\bar X_{t_k}^\pi|^2 \Big]  & \leq & 2   \E\Big[\esssup_{a\in A} \E_{t_1,a}\big[ \ldots \esssup_{a\in A}  
\E_{t_k,a}\big[ \sup_{t\in[t_k,t_{k+1}]}\big| X_t  -  X_{t_k} \big|^2 \big] \ldots \big] \Big] \\
& & \; + \; 2 \E\Big[\esssup_{a\in A} \E_{t_1,a}\big[ \ldots \esssup_{a\in A}   \E_{t_{k-1},a}\big[ |\Delta X_{t_{k}}|^2  \big] \ldots \big] \Big] \\
& \leq & C_1 |\pi|,
\enqs
by \reff{interXpi}-\reff{interXpi2}, which ends the proof.  
\ep

\section{Discretely jump-constrained BSDE}

\setcounter{equation}{0} \setcounter{Assumption}{0}
\setcounter{Theorem}{0} \setcounter{Proposition}{0}
\setcounter{Corollary}{0} \setcounter{Lemma}{0}
\setcounter{Definition}{0} \setcounter{Remark}{0}

Given the  forward regime switching process $(X,I)$ defined in the previous section, we  consider the minimal quadruple solution  $(Y,Z,U,K)$  to the  BSDE with
nonpo\-si\-ti\-ve jumps:
\begin{equation} \label{BSDEnonpos2}
\left\{ \begin{array}{ccl}
Y_t &=& g(X_T) +  \int_t^T f(X_s,I_s,Y_s,Z_s) ds + K_T - K_t \\
& & \;\; - \int_t^T Z_s dW_s - \int_t^T \int_A U_s(a) \tilde\mu(ds,da),  \;\;\;\;\;   0 \leq t \leq T.  \\
U_t(a) & \leq & 0,  
\end{array}
\right.
\end{equation}
By  solution to \reff{BSDEnonpos2}, we mean a quadruple $(Y,Z,U,K)$  $\in$ $\Sc^2\times L^2(W)\times L^2(\tilde\mu)\times\Kc^2$, where  $\Sc^2$ is  the space of c\`ad-l\`ag or c\`ag-l\`ad $\F$-progressively measurable processes $Y$ satisfying $\|Y\|^2$ $:=$ 
$\E[\sup_{t\in[0,T]}|Y_t|^2]$ $<$ $\infty$,   $L^2(W)$ is   the space of $\R^d$-valued $\Pc$-measurable processes  such that $\|Z\|^2_{L^2(W)}$ $:=$ $\E[\int_0^T|Z_t|^2dt]$ $<$ $\infty$,  
$L^2(\tilde\mu)$ is  the space of  real-valued $\Pc\otimes \Bc(A)$-measurable   processes $U$ such that  $\|U\|^2_{L^2(\tilde\mu)}$ $:=$ $\E[\int_0^T\int_A|U_t(a)|^2\lambda(da)\,dt]$ $<$ $\infty$,  and $\Kc^2$ is   the subspace of $\Sc^2$ consisting of nondecreasing predictable processes  such that $K_0=0$, $\P$-a.s., and  the equation in  \eqref{BSDEnonpos2} holds $\P$-a.s., while  the nonpositive jump constraint holds on 
$\Omega\times [0,T]\times A$ a.e. with respect to the measure $d\P\otimes dt\otimes \lambda(da)$. 
By minimal solution to the BSDE \eqref{BSDEnonpos}, we mean a quadruple solution $(Y,Z,U,K)$ $\in$ $\Sc^2\times L^2(W)\times L^2(\tilde\mu)\times\Kc^2$  such that for any other solution
$(Y',Z',U',K')$ to the same BSDE, we have $\P$-a.s.: $Y_t$  $\leq$ $Y_t'$,  $t \in [0,T]$.  

\vspace{1mm}

In the rest of this paper, we shall make the standing Lipschitz assumption on the functions $f$ $:$ $\R^d\times A\times\R\times\R^d$ $\to$ $\R$ and  $g$ $:$ $\R^d$ $\to$ $\R$.

\vspace{2mm}

\ni  \textbf{(H2)}
There exists a constant $L_2$ such that
\beqs
 |f(x,a,y,z) - f(x',a',y',z')|+ |g(x) - g(x')| & \leq &  L_2 \big(|x-x'| +|a-a'| +|y-y'| + |z-z'|\big), 
\enqs
for all $x,x'\in \R^d$, $y,y' \in \R$, $z,z' \in \R^d$, $a,a'\in A$.  In the sequel, we shall denote by $C$ a generic positive constant  which depends only on $L_1$, $L_2$, $T$, $(X_0,I_0)$ 
and  $\lambda(A)$ $<$ $\infty$, and may vary from lines to lines.

\vspace{2mm}

Under {\bf (H1)}-{\bf (H2)}, it is proved in \cite{kp12} the existence and uniqueness of a minimal solution $(Y,Z,U,K)$ to \reff{BSDEnonpos2}.  Moreover, the minimal solution $Y$ is in the form 
\beq \label{Yv} 
Y_t &=& v(t,X_t), \;\;\; 0 \leq t \leq T,
\enq
where $v$ $:$ $[0,T]\times\R^d$ $\to$ $\R$  is a viscosity solution with linear growth to the fully nonlinear HJB type equation: 
\begin{equation} \label{HJB} 
\left\{ \begin{array}{rcl} 
-  \sup_{a \in A} \big [ \Lc^a v + f(x,a,v,\sigma\trans(x,a)D_x v) \big] &=& 0,   \mbox{ on } [0,T) \times \R^d, \\
v(T,x) &=& g,   \mbox{ on } \R^d,
\end{array}
\right.
\end{equation}
where
\beqs
\Lc^a v &=& \Dt{v} + b(x,a). D_x v + \frac{1}{2} {\rm tr}(\sigma\sigma\trans(x,a)D_x^2 v).  
\enqs

We shall make the standing assumption  that comparison principle holds for \reff{HJB}.  

\vspace{2mm}

\ni {\bf (HC)} Let $\bar w$ (resp. $\underline{w}$) be a lower-semicontinuous (resp. upper-semicontinuous) viscosity supersolution (resp. subsolution) with linear growth condition  to \reff{HJB}. Then, 
$\bar w$ $\geq$ $\underline{w}$.

\vspace{1mm}

When $f$ does not depend on $y,z$, i.e. \reff{HJB} is the usual HJB equation for a stochastic control problem,  Assumption {\bf (HC)} holds true, see \cite{FS} or \cite{pha09}.  In the general case, we refer 
to \cite{craishlio92} for sufficient conditions to comparison principles. Under {\bf (HC)}, the function $v$ in \reff{Yv}  is the unique viscosity solution to \reff{HJB}, and is in particular continuous. Actually, we have 
the standard H\"older and Lipschitz property (see Appendix in \cite{kr00} or \cite{barjac07}):
\beq \label{vlip}
|v(t,x) - v(t',x')| & \leq & C \big( |t-t'|^{1\over 2} + | x- x'|\big), \;\;\;  (t,t') \in [0,T], x,x' \in \R^d.
\enq
This implies  that the process $Y$ is continuous, and thus  the jump component  $U$ $=$ $0$.  
In the sequel, we shall focus on  the approximation of the remaining components  $Y$ and $Z$ of the minimal solution to \reff{BSDEnonpos2}.

\vspace{1mm}

\vspace{2mm}

We introduce in this section  discretely jump-constrained BSDE.   The nonpositive jump constraint operates only at the times of the grid $\pi=\{t_0=0<t_1<\ldots<t_n=T\}$ of $[0,T]$, and  
we look for a quadruple $(Y ^\pi,  \Yc ^\pi, \Zc ^\pi, \Uc ^\pi)$ $\in$ 
$\Sc^2\times\Sc^2\times L^2(W)\times L^2(\tilde\mu)$ satisfying:
\beq\label{BSDEDC1}
 Y^\pi_T & = &  \Yc^\pi_T~=~g(X_T)
 \enq
 and 
 \beq\label{BSDEDC2}
 \Yc ^\pi_t & = &  Y ^\pi_{t_{k+1}} + \int_t^{t_{k+1}}f(X_s,I_s, \Yc ^\pi_s, \Zc ^\pi_s)ds\\\nonumber
 & & -\int_t^{t_{k+1}} \Zc ^\pi_sdW_s -\int_t^{t_{k+1}}\int_A  \Uc ^\pi_s(a)\tilde \mu(ds,da)\;,
 \\\label{BSDEDC3}
 Y ^\pi_{t} & = & \Yc ^\pi_{t}\mathds{1}_{(t_k,t_{k+1})}(t)+ \esssup_{a\in A} \E\Big[ \Yc ^\pi_{t}\big| X_{t},I_{t}~=~a\Big]\mathds{1}_{\{t_k\}}(t)\;,\qquad 
\enq
for all $t\in[t_k,t_{k+1})$ and all $ 0\leq k\leq n -1$.

\vspace{2mm}

Notice that at each time $t_k$ of the grid, the condition is not known a priori to be square integrable since it involves a supremum over $A$, and the well-posedness of the BSDE \reff{BSDEDC1}-\reff{BSDEDC2}-\reff{BSDEDC3} is not a direct and standard issue. We shall use a PDE approach for proving the existence and uniqueness of a solution.  Let us consider the system of integro-partial differential equations (IPDEs) for the functions $v^\pi$ and $\vartheta^\pi$  defined recursively on $[0,T]\times\R^d\times A$ by:

 \noindent $\bullet$ A  terminal condition for $v^\pi$ and $\vartheta^\pi$:
\beq\label{Term-cond}
 v^\pi (T,x,a)~~=~~ \vartheta^\pi (T,x,a) & = & g(x) \;,\quad (x,a)\in \R^d\times A\;,
\enq

\noindent $\bullet$ A sequence of IPDEs for $\vartheta^\pi$ 
\begin{equation}\left\{
\begin{array}{rcl} 
  -  \Lc^a \vartheta ^\pi - f\big(x,a,\vartheta ^\pi,\sigma\trans(x,a)D_x\vartheta ^\pi\big) & &  \\
- \int_A\big(\vartheta ^\pi(t,x,a')-\vartheta ^\pi(t,x,a)\big)\lambda(da') \; = \;   0, & &   (t,x,a)\in [t_k,t_{k+1})\times\R^d\times A, \\
\vartheta^\pi (t_{k+1}^-,x,a) \; = \;  \sup_{a'\in A}\vartheta ^\pi (t_{k+1},x,a') & &  (x,a)\in \R^d\times A 
\end{array}\right.
\label{seqIPDE}  
\end{equation}
for $k=0\ldots,n -1$,

\vspace{1mm}

\noindent $\bullet$ the relation between $v^\pi$ and $\vartheta^\pi$:
\beq \label{relvvartheta}
v^\pi(t,x,a) & = & \vartheta ^\pi({t},x,a)\mathds{1}_{(t_k,t_{k+1})}(t)+ \sup_{a'\in A} \vartheta ^\pi({t},x,a')\mathds{1}_{\{t_k\}}(t)\;,\qquad 
\enq
for all $t\in[t_k,t_{k+1})$ and $k=0\ldots,n -1$.  
The rest of this section is devoted to the proof of existence and uniqueness of a solution to \reff{Term-cond}-\reff{seqIPDE}-\reff{relvvartheta}, together with some uniform Lipschitz properties, and its connection to the discretely jump-constrained BSDE \reff{BSDEDC1}-\reff{BSDEDC2}-\reff{BSDEDC3}.

 \vspace{2mm}

For any  $L$-Lipschitz continuous function $\varphi$ on $\R^d\times A$,   and $k$ $\leq$ $n-1$,  we denote: 
\beq \label{operaT}
\T_\pi^k[\varphi](t,x,a) & := & w(t,x,a), \;\;\; (t,x,a) \in [t_k,t_{k+1})\times\R^d\times A, 
\enq
where $w$ is the unique continuous viscosity solution on $[t_k,t_{k+1}]\times\R^d\times A$  with linear growth condition in $x$ to the integro partial differential equation (IPDE): 
\begin{equation}\left\{
\begin{array}{rcl}  
 -  \Lc^a  w  - f(x,a,w,\sigma\trans D_x w) & &  \\
-\int_A\big(w(t,x,a')-w(t,x,a)\big)\lambda(da') \; = \;  0,  & &  (t,x,a)\in [t_k,t_{k+1})\times\R^d\times A, \\
w (t_{k+1}^-,x,a) \; = \; \varphi(x,a), & &  (x,a)  \in \R^d\times A\;, \label{IPDEbarn}
\end{array}\right.
\end{equation}
and we extend by continuity $\T_\pi^k[\varphi](t_{k+1},x,a)$ $=$ $\varphi(x,a)$. 
The existence and uniqueness of such a solution $w$ to the semi linear IPDE \reff{IPDEbarn}, and its nonlinear Feynman-Kac representation in terms of   
BSDE with jumps, is obtained e.g. from  Theorems 3.4 and 3.5 in \cite{bbp97}.

\vspace{2mm}

\begin{Lemma}\label{lem stab lip}
There exists a constant $C$ such that for any $L$-Lipschitz continuous function $\varphi$ on $\R^d\times A$,  and $k$ $\leq$ $n-1$, 
we have
\beqs
|\T_\pi^k[\varphi](t,x,a)-\T_\pi^k[\varphi](t,x',a')| & \leq & \max(L,1) \sqrt{1+|\pi|}e^{C|\pi|}(|x-x'|+|a-a'|)\;,
\enqs 
for all  $t\in[t_k,t_{k+1})$, and $(x,a),(x',a')\in\R^d\times A$.
\end{Lemma}
\textbf{Proof.} Fix $t\in [t_k,t_{k+1})$, $k$ $\leq$ $n-1$,  $(x,a), (x',a')\in \R^d\times A$, and $\varphi$ an $L$-Lipschitz continuous function on $\R^d\times A$. 
Let $(Y^\varphi,Z^\varphi,U^\varphi)$ and $(Y^{\varphi,'},Z^{\varphi,'},U^{\varphi,')}$ be the solutions on $[t,t_{k+1}]$ 
to the BSDEs
\beqs
Y_s^\varphi & = & \varphi(X^{t,x,a}_{t_{k+1}},I_{t_{k+1}}^{t,a})+\int_s^{t_{k+1}}f(X_r^{t,x,a},I^{t,a}_r,Y_r^\varphi,Z_r^\varphi)dr \\
& & \;\;\;  - \; \int_s^{t_{k+1}}Z_r^\varphi dW_r-\int_s^{t_{k+1}}\int_AU_r^\varphi(e)\tilde \mu (dr,de),  \;\;\;\;\;  t \leq s \leq t_{k+1}, \\
Y^{\varphi,'}_s & = & \varphi(X^{t,x',a'}_{t_{k+1}},I_{t_{k+1}}^{t,a'}) + \int_s^{t_{k+1}}f(X^{t,x',a'}_r,I^{t,a'}_r,Y^{\varphi,'}_r,Z^{\varphi,'}_r)dr \\
& & \;\;\; - \; \int_s^{t_{k+1}}Z^{\varphi,'}_r dW_r-\int_s^{t_{k+1}}\int_AU^{\varphi,'}_r(e)\tilde \mu (dr,de),  \;\;\;\;\;  t \leq s \leq t_{k+1}
\enqs
From Theorems 3.4 and 3.5 in \cite{bbp97}, we have the identification: 
\beq\label{ide bsde tnpsi}
Y_t^\varphi~~=~~\T_\pi^k[\varphi](t,x,a)  & \mbox{ and } & Y^{\varphi,'}_t~~=~~\T_\pi^k[\varphi](t,x',a')\;.
\enq
We now estimate the difference between the processes $Y^\varphi$and $Y^{\varphi,'}$, and  set 
$\delta Y^\varphi$ $=$ $Y^\varphi-Y^{\varphi,'}$, $\delta Z^\varphi$ $=$ $Z^\varphi-Z^{\varphi,'}$, $\delta X$ $=$ 
$X^{t,x,a}-X^{t,x',a'}$, $\delta I$ $=$ $I^{t,a}-I^{t,a'}$. 
By  It\^o's formula, the Lipschitz condition of $f$ and $\varphi$, and Young inequality,  we have
\beqs
\E\Big[|\delta Y_s^\varphi|^2\Big]  + \E \Big[ \int_s^{t_{k+1}} |\delta Z_s^\varphi|^2 ds \Big] 
& \leq & L^2  \E\Big[|\delta X_T|^2 + |\delta I_T|^2\Big] +  C\int_s^{t_{k+1}}\E\Big[|\delta Y_r^\varphi|^2\Big]dr\\
 & & + \frac{1}{2} \E\Big[ \int_s^{t_{k+1}} \big( |\delta X_r|^2+|\delta I_r|^2 + | \delta Z_r^\varphi|^2 \big) dr  \Big],
\enqs
for any $s\in [t,t_{k+1}]$. Now, from classical estimates on jump-diffusion processes we have
\beqs
\sup_{r\in[t,t_{k+1}]}\E\Big[|\delta X_r|^2 + |\delta I_r|^2\Big] & \leq & e^{C|\pi|}\big(|x-x'|^2+|a-a'|^2\big), 
\enqs
and thus:
\beqs
\E\Big[|\delta Y_s^\varphi|^2\Big] &\leq  & (L^2+|\pi|)e^{C|\pi|}\big(|x-x'|^2+|a-a'|^2\big) + C\int_s^{t_{k+1}}\E\Big[|\delta Y_r^\varphi|^2\Big]dr\;,
\enqs
for all $s\in [t,t_{k+1}]$. By Gronwall's Lemma, this yields
\beqs
\sup_{s\in[t,t_{k+1}]}\E\Big[|\delta Y_s^\varphi|^2\Big] &\leq & (L^2+|\pi|)e^{2C|\pi|}\big(|x-x'|^2+|a-a'|^2\big),
\enqs
which proves the required result from the identification \reff{ide bsde tnpsi}: 
\beqs
|\T_\pi^k[\varphi](t,x,a)-\T_\pi^k[\varphi](t,x',a')| & \leq & \sqrt{L^2+|\pi|}e^{C|\pi|}(|x-x'|+|a-a'|)\\
 & \leq & \max(L,1)\sqrt{1+|\pi|}e^{C|\pi|}(|x-x'|+|a-a'|).
\enqs
\ep

\vspace{2mm}


\begin{Proposition} 
There exists a unique viscosity solution $\vartheta^\pi$ with linear growth condition to the  IPDE \reff{Term-cond}-\reff{seqIPDE}, and this solution satisfies:
\beq
& & |\vartheta^\pi(t,x,a)-\vartheta^\pi(t,x',a')|  \nonumber \\
& \leq &  \max(L_2,1) \sqrt{\Big(e^{2C|\pi|}(1+|\pi| )\Big)^{n -{k}}}\big(|x-x'|+|a-a'|\big),  \label{ppte lipschitz vartheta pi}
\enq
for all $k$ $=$ $0,\ldots,n-1$, $t\in [t_k,t_{k+1})$,  $(x,a),(x',a')\in \R^d\times A$.
\end{Proposition}
\textbf{Proof.} We prove by a backward induction on $k$ that the IPDE \reff{Term-cond}-\reff{seqIPDE} admits a unique solution on 
$[t_k,T]\times \R^d\times A$, which satisfies \reff{ppte lipschitz vartheta pi}. 

\vspace{1mm}

\ni $\bullet$ For $k=n-1$,  we directly get the existence and uniqueness of $\vartheta^\pi$ on $[t_{n-1},T]\times\R^d\times A$ from  Theorems 3.4 and 3.5 in \cite{bbp97}, and we have $\vartheta^\pi$ $=$ $\T_\pi^{n-1}[g]$ on $[t_{n-1},T)\times\R^d\times A$. Moreover, we also get by Lemma  
\ref{lem stab lip}:
\beqs
|\vartheta^\pi(t,x,a)-\vartheta^\pi(t,x',a')| & \leq &  \max(L_2,1) \sqrt{e^{2C|\pi|}(1+|\pi| )}\big(|x-x'|+|a-a'|\big)
\enqs
for all $t\in [t_{n-1},t_n)$, $(x,a),(x',a')\in \R^d\times A$.

\vspace{1mm}

\ni $\bullet$ Suppose that the result holds true at step  $k+1$ i.e. there exists a unique function $\vartheta^\pi$ on $[t_{k+1},T]\times\R^d\times A$ with linear growth and satisfying \reff{Term-cond}-\reff{seqIPDE} and \reff{ppte lipschitz vartheta pi}.   It remains to prove that $\vartheta^\pi$ is uniquely determined  by \reff{seqIPDE} on $[t_k,t_{k+1})\times \R^d\times A$ and that it satisfies \reff{ppte lipschitz vartheta pi} on 
$[t_k,t_{k+1})\times \R^d\times A$. Since $\vartheta^\pi$ satisfies \reff{ppte lipschitz vartheta pi} at time $t_{k+1}$, we deduce that the function 
\beqs
\psi_{k+1}(x) &:=&  \sup_{a\in A} \vartheta^\pi(t_{k+1},x,a), \;\;\; x \in \R^d,
\enqs
is also Lipschitz continuous, and satisfies by the induction hypothesis:
\beq \label{psilip}
|\psi_{k+1}(x)-\psi_{k+1}(x')| & \leq & \max(L_2,1)  \sqrt{\Big(e^{2C|\pi|}(1+|\pi| )\Big)^{n -{k-1}}}|x-x'|, 
\enq
for all $x,x'\in\R^d$.  Under {\bf (H1)} and {\bf (H2)}, we can apply Theorems 3.4 and 3.5 in \cite{bbp97}, and we get that $\vartheta^\pi$ is the unique 
viscosity solution with linear growth to \reff{seqIPDE} on $[t_k,t_{k+1})\times \R^d\times A$, with $\vartheta^\pi$ $=$ $\T_\pi^k[\psi_{k+1}]$. 
Thus it exists and is unique on  $[t_k,T]\times \R^d\times A$.   From Lemma \ref{lem stab lip} and \reff{psilip},  we then get
\beqs
|\vartheta^{\pi}(t,x,a)-\vartheta^{\pi}(t,x',a')| & = & |\T_\pi^k[\psi_{k+1}](t,x,a)-\T_\pi^k[\psi_{k+1}](t,x',a')|\\
 & \leq & \max(L_2,1) \sqrt{\Big(e^{2C|\pi|}(1+|\pi| )\Big)^{n -{k-1}}}  \\
 & &  \;  . \sqrt{(1+|\pi|)e^{2C|\pi|}}| \big(|x-x'|+|a-a'|\big)\\
 & \leq & \max(L_2,1)  \sqrt{\Big(e^{2C|\pi|}(1+|\pi| )\Big)^{n -{k}}} \big(|x-x'|+|a-a'|\big)
\enqs
for any $t\in [t_k,t_{k+1})$ and $(x,a),(x',a')\in \R^d\times A$, which  proves the required induction inequality at step $k$. 
\ep

\begin{Remark} \label{remvarthetapi}
{\rm  The function $a$ $\to$ $\vartheta^\pi(t,x,.)$ is continuous on $A$, for each $(t,x)$, and so the function $v^\pi$ is  well-defined by \reff{relvvartheta}.  
Moreover, the function $\vartheta^\pi$ may be written recursively as: 
\begin{equation} \label{varthetapi}
\left\{ 
\begin{array}{rcll}
\vartheta^\pi(T,.,.) &=&  g &  \mbox{ on } \;  \R^d\times A, \\
\vartheta^\pi &=& \T_\pi^k[v^\pi(t_{k+1},.)],  & \mbox{ on } \; [t_{k},t_{k+1})\times\R^d\times A, 
\end{array}
\right.
\end{equation}
for $k = 0,\ldots,n-1$.  In particular, $\vartheta^\pi$  is continuous on $(t_k,t_{k+1})\times\R^d\times A$, $k$ $\leq$ $n-1$.    
\ep
}
\end{Remark}

\vspace{3mm}

As a consequence of the above proposition,  we obtain the uniform Lipschitz property of $\vartheta^\pi$ and $v^\pi$, with a Lipschitz constant independent of $\pi$.

\begin{Corollary} \label{prop stab lip}
There exists a constant $C$ (independent of $|\pi|$)  such that
\beqs
|\vartheta^\pi(t,x,a)-\vartheta^\pi(t,x',a')|  +  |v^\pi(t,x,a)-v^\pi(t,x',a')|  & \leq & C \big(|x-x'|+|a-a'|\big),
\enqs
for all  $t\in [0,T]$, $x,x'\in \R^d$, $a,a'\in \R^d$. 
\end{Corollary}
\textbf{Proof.} Recalling that $n|\pi|$ is bounded, we see that  the sequence appearing in \reff{ppte lipschitz vartheta pi}: 
 $\Big(\big(e^{2C|\pi|}(1+|\pi| )\big)^{n -{k}}\Big)_{0\leq k\leq n -1}$ is bounded uniformly in $|\pi|$ (or $n$), 
 which shows the required Lipschitz property of $\vartheta^\pi$.  Since $A$ is assumed to be compact, this shows in particular that the function 
 $v^\pi$ defined by the relation  \reff{relvvartheta} is  well-defined and finite. Moreover, by noting that 
\beqs
|\sup_{a\in A}\vartheta^\pi(t,x,a)-\sup_{a\in A}\vartheta^\pi(t,x',a)| & \leq  & \sup_{a\in A}|\vartheta^\pi(t,x,a)-\vartheta^\pi(t,x',a)|
\enqs
for all $(t,x)\in[0,T]\times\R^d$, we also obtain  the required Lipschitz property for $v^\pi$. 
 \ep

 \vspace{3mm}

We now turn to the existence of a solution to the discretely jump-constrained BSDE.

\begin{Proposition}\label{exist uniq EDSR disc const}
The BSDE \reff{BSDEDC1}-\reff{BSDEDC2}-\reff{BSDEDC3} admits a unique solution $(Y^\pi,\Yc^\pi,\Zc^\pi,\Uc^\pi)$ in $\Sc^2\times \Sc^2\times L^2(W)\times L^2(\tilde\mu)$. Moreover we have
\beq\label{identif BSDE IPDE}
\Yc^\pi_t ~~ = ~~ \vartheta^\pi(t,X_t,I_t), & \mbox{ and } & Y^\pi_t ~~ = ~~ v^\pi(t,X_t,I_t) 
\enq
for all $t\in[0,T]$.
\end{Proposition}
\textbf{Proof.} We  prove by backward induction on $k$ that $(Y^\pi,\Yc^\pi,\Zc^\pi,\Uc^\pi)$ is well defined and satisfies \reff{identif BSDE IPDE} on 
$[t_k,T]$. 

\vspace{1mm}

\ni $\bullet$ Suppose that $k=n-1$. From Corollary 2.3 in \cite{bbp97}, we know that $(\Yc^\pi,\Zc^\pi,\Uc^\pi)$,  exists and is unique on  $[t_{n-1},T]$. Moreover,  from Theorems 3.4 and 3.5 in \cite{bbp97}, we get $\Yc_t^\pi$ $=$ $\T_\pi^k[g](t,X_t,I_t)$ $=$ $\vartheta^\pi(t,X_t,I_t)$ on $[t_{n-1},T]$.  By \reff{BSDEDC3}, we then have for all $t$ $\in$ $[t_{n-1},T)$:
\beqs
Y_t^\pi &=& \mathds{1}_{(t_{n-1},T)}(t) \; \vartheta^\pi(t,X_t,I_t)  +  \mathds{1}_{t_{n-1}}(t) \; \esssup_{a\in A} \vartheta^\pi(t,X_t,a)  \\
&=& \mathds{1}_{(t_{n-1},T)}(t) \; \vartheta^\pi(t,X_t,I_t)  +  \mathds{1}_{t_{n-1}}(t) \; \sup_{a\in A} \vartheta^\pi(t,X_t,a)  \; = \; v^\pi(t,X_t,I_t),
\enqs
since the essential supremum and supremum coincide by continuity of $a$ $\to$ $\vartheta^\pi(t,X_t,a)$ on the compact set $A$. 

\vspace{1mm}

\ni $\bullet$ Suppose that the result holds true for some $k\leq n-1$. Then, we see that  $(\Yc^\pi,\Zc^\pi,\Uc^\pi)$ is defined on $[t_{k-1},t_k)$ as the solution to a BSDE driven by $W$ and 
$\tilde \mu$ with a terminal condition 
$v^\pi(t_k,X_{t_k})$. Since $v^\pi$ satisfies a linear growth condition,  
we know again by Corollary 2.3 in \cite{bbp97}  that $(\Yc^\pi,\Zc^\pi,\Uc^\pi)$, thus also $Y^\pi$,  exists and is unique on  $[t_{k-1},t_k)$. 
Moreover, using again Theorems  3.4 and 3.5 in \cite{bbp97},  we get \reff{identif BSDE IPDE} on $[t_{k-1},t_k)$. 
\ep



\vspace{2mm}

We end this section with a conditional regularity result for the discretely jump-constrained BSDE.

\begin{Proposition} \label{regulYdis}
There exists some constant $C$ such that
\beqs
 \sup_{t \in [t_k,t_{k+1})} \E_{t_k}\big[ |\Yc_t^\pi - \Yc_{t_k}^\pi|^2 \big]  +  
 \sup_{t \in (t_k,t_{k+1}]} \E_{t_k}\big[ |Y_t^\pi - Y_{t_{k+1}}^\pi|^2 \big]  
 & \leq & C(1+ |X_{t_k}|^2) |\pi|,
\enqs
for all $k$ $=$ $0,\ldots,n-1$. 
\end{Proposition}
{\bf Proof.} Fix $k$ $\leq$ $n-1$. By It\^o's formula, we have for all $t$ $\in$ $[t_k,t_{k+1})$:
\beqs
\E_{t_k}\big[ |\Yc_t^\pi - \Yc_{t_k}^\pi|^2 \big] &=& 2 \E_{t_k}\Big[ \int_{t_k}^t  f(X_s,I_s,\Yc_s^\pi,\Zc_s^\pi) (\Yc_{t_k}^\pi - \Yc_s^\pi) ds \Big]  \\
& & \;\;\; + \; \E_{t_k}\Big[ \int_{t_k}^t |\Zc_s^\pi|^2 \Big] + \E_{t_k}\Big[ \int_{t_k}^t \int_A |\Uc_s^\pi(a)|^2\lambda(da) ds \Big] \\
& \leq & \E_{t_k}\Big[ \int_{t_k}^t |\Yc_s^\pi - \Yc_{t_k}^\pi|^2 \Big]  + C |\pi| \Big( 1 +  \E_{t_k}\Big[ \sup_{s\in[t_k,t_{k+1}]}|X_s|^2 \Big] \Big) \\
& & + \;  C |\pi|  \E_{t_k}\Big[ \sup_{s\in[t_k,t_{k+1}]}  \Big( |\Yc_s^\pi|^2 + |\Zc_s^\pi|^2 + \int_A |\Uc_s^\pi(a)|^2\lambda(da) \Big) \Big],
\enqs
by the linear growth condition on $f$ (recall also that $A$ is compact), and Young inequality. 
Now, by standard estimate  for $X$ under growth linear condition on $b$ and $\sigma$, we have:
\beq \label{estimXtk}
 \E_{t_k}\Big[ \sup_{s\in[t_k,t_{k+1}]}|X_s|^2 \Big]  & \leq & C(1 + |X_{t_k}|^2). 
\enq
We also know from Proposition 4.2 in \cite{BE08}, under {\bf (H1)} and {\bf (H2)}, that  there exists a constant $C$ depending only on the Lipschitz constants of $b$, $\sigma$ $f$ and $v^\pi(t_{k+1},.)$ (which does not depend on $\pi$ by Corollary \ref{prop stab lip}), such that
\beq \label{estimYZU}
\E_{t_k}\Big[ \sup_{s\in[t_k,t_{k+1}]} \Big( |\Yc_s^\pi|^2 + |\Zc_s^\pi|^2 + \int_A |\Uc_s^\pi(a)|^2\lambda(da) \Big) \Big] & \leq &   C(1+|X_{t_k}|^2). 
\enq
We deduce that
\beqs
\E_{t_k}\big[ |\Yc_t^\pi - \Yc_{t_k}^\pi|^2 \big]  & \leq &  \E_{t_k}\Big[ \int_{t_k}^t |\Yc_s^\pi - \Yc_{t_k}^\pi|^2 \Big]  + C |\pi| (1 + |X_{t_k}|^2),
\enqs
and we conclude  for the regularity of $\Yc^\pi$ by Gronwall's lemma. Finally, from the definition \reff{BSDEDC2}-\reff{BSDEDC3} 
of $Y^\pi$ and $\Yc^\pi$, It\^o isometry for stochastic integrals, and growth linear condition on $f$, we have for all $t$ $\in$ $(t_k,t_{k+1})$:
\beqs
\E_{t_k}\big[ |Y_t^\pi - Y_{t_{k+1}}^\pi|^2 \big] &=&  \E_{t_k}\big[ |\Yc_t^\pi - Y_{t_{k+1}}^\pi|^2 \big]  \\
& \leq & 3 \E_{t_k} \Big[ \int_{t_k}^{t_{k+1}} \Big( |f(X_s,I_s,\Yc_s^\pi,\Zc_s^\pi)|^2  + |\Zc_s^\pi|^2 +  \int_A |\Uc_s^\pi(a)|^2\lambda(da) \Big) ds \Big] \\
& \leq & C |\pi|  \E_{t_k} \Big[ 1 +  \sup_{s\in[t_k,t_{k+1}]}  \Big(  |X_s|^2 + |\Yc_s^\pi|^2 + |\Zc_s^\pi|^2 + \int_A |\Uc_s^\pi(a)|^2\lambda(da) \Big) \Big] \\
& \leq & C |\pi| (1 + |X_{t_k}|^2),
\enqs
where we used again \reff{estimXtk} and \reff{estimYZU}. This ends the proof. 
\ep

\section{Convergence of discretely jump-constrained BSDE}

\setcounter{equation}{0} \setcounter{Assumption}{0}
\setcounter{Theorem}{0} \setcounter{Proposition}{0}
\setcounter{Corollary}{0} \setcounter{Lemma}{0}
\setcounter{Definition}{0} \setcounter{Remark}{0}

This section is devoted to the convergence of the discretely jump-constrained BSDE towards the minimal solution to the BSDE with nonpositive jump.

\subsection{Convergence result}

\begin{Lemma} \label{propcroisbarvn}
We have the following assertions:

\ni 1) The familly $(\vartheta^\pi)_\pi$ is nondecreasing and upper bounded by $v$: 
for any grids $\pi$ and $\pi'$ such that $\pi$ $\subset$ $\pi'$, we have  
\beqs
\vartheta^\pi(t,x,a) & \leq & \vartheta^{\pi'}(t,x,a)~\leq~v(t,x)\;,\quad (t,x,a)\in [0,T]\times \R^d\times A\;. 
\enqs
2) The familly $(\vartheta^\pi)_\pi$ satisfies a uniform linear growth condition: there exists a constant $C$ such that
\beqs
|\vartheta^\pi(t,x,a)| & \leq & C(1+|x|), 
\enqs
for any $(t,x,a)\in[0,T]\times\R^d\times A$ and any grid $\pi$.
\end{Lemma}
\textbf{Proof.}  1)  Let us first prove that $\vartheta^\pi$ $\leq$ $v$.  Since $v$ is a (continuous) viscosity solution to the HJB equation \reff{HJB}, 
and $v$ does not depend on $a$,  we see that $v$ is a viscosity supersolution 
to the IPDE in \reff{seqIPDE} satisfied by $\vartheta^\pi$ on each interval $[t_k,t_{k+1})$. Now, since $v(T,x)$ $=$ $\vartheta^\pi(T,x,a)$, we deduce  by comparison  principle for this IPDE 
(see e.g. Theorem 3.4 in \cite{bbp97}) on  $[t_{n-1},T)\times\R^d\times A$  that  $v(t,x)$ $\geq$ $\vartheta^\pi(t,x,a)$ for all $t$ $\in$ $[t_{n-1},T]$, $(x,a)$ $\in$ $\R^d\times A$. In particular,  $v(t_{n-1}^-,x)$ $=$ $v(t_{n-1},x)$ $\geq$ $\sup_{a\in A}\vartheta^\pi(t_{n-1},x,a)$ $=$  $\vartheta^\pi(t_{n-1}^-,x,a)$. Again, by comparison principle for the IPDE \reff{seqIPDE} on $[t_{n-2},t_{n-1})\times\R^d\times A$, it follows that  
$v(t,x)$ $\geq$ $\vartheta^\pi(t,x,a)$ for all $t$ $\in$ $[t_{n-2},t_{n-1}]$, $(x,a)$ $\in$ $\R^d\times A$.  By backward induction on time, we conclude  that $v$ $\geq$ $\vartheta^\pi$ on $[0,T]\times\R^d\times A$.  

Let us next consider two partitions  $\pi$ $=$ $(t_k)_{0\leq k\leq n}$ and $\pi'$ $=$ $(t'_k)_{0\leq k\leq n'}$ of $[0,T]$ with $\pi$ $\subset$ $\pi'$, and denote by $m$ $=$ $\max\{ k \leq n': t_{m}' \notin \pi\}$. 
Thus, all the points of the grid $\pi$ and $\pi'$ coincide after time $t'_m$, and  since $\vartheta^\pi$ and $\vartheta^{\pi'}$  are viscosity solution to the same IPDE  \reff{seqIPDE} 
starting from the same terminal data $g$, we deduce 
by uniqueness that $\vartheta^\pi$ $=$ $\vartheta^{\pi'}$ on $[t'_m,T]\times\R^d\times A$.  Then, we have $\vartheta^{\pi'}(t_m^{'-},x,a)$ $=$ $\sup_{a\in A}\vartheta^\pi(t_m',x,a)$ $=$ 
$\sup_{a \in A} \vartheta^\pi(t_m',x,a)$ $\geq$ $\vartheta^\pi(t_m^-,x,a)$ since $\vartheta^\pi$ is continuous outside of the points of the grid $\pi$ (recall Remark \ref{remvarthetapi}).  Now, 
since $\vartheta^\pi$ and $\vartheta^{\pi'}$  are viscosity solution to the same IPDE \reff{seqIPDE} on $[t'_{m-1},t_m)$, we deduce by comparison principle 
that $\vartheta^{\pi'}$ $\geq$ $\vartheta^{\pi}$ on $[t'_{m-1},t_m']\times\R^d\times A$.  Proceeding by backward induction, we conclude that $\vartheta^{\pi'}$ $\geq$ $\vartheta^{\pi}$ on $[0,T]\times\R^d\times A$.

\vspace{1mm} 

\ni 2)  Denote by $\pi_0=\{t_0=0,t_1=T\}$ the trivial grid of $[0,T]$. Since $\vartheta^{\pi_0}\leq \vartheta^{\pi}\leq v$  and $\vartheta^{\pi_0}$ and $v$ satisfy a linear growth condition, we get (recall that $A$ is compact):
\beqs
|\vartheta^{\pi}(t,x,a)| & \leq & |\vartheta^{\pi_0}(t,x,a)|+|v(t,x)| \; \leq \;  C(1+|x|),
\enqs 
for any $(t,x,a)\in[0,T]\times\R^d\times A$ and any grid $\pi$. 
\ep

\vspace{3mm}

In the sequel, we denote by $\vartheta$ the increasing limit of the sequence $(\vartheta^\pi)_\pi$ when the grid increases by becoming finer, i.e. its modulus  $|\pi|$ goes to zero.  The next result shows that $\vartheta$ does not depend on the variable $a$ in $A$.

\begin{Proposition}\label{prop limv indepa} 
The function $\vartheta$ is l.s.c. and does not depend on the variable $a\in A$: 
\beqs
\vartheta (t,x,a) & = & \vartheta(t,x,a')\;,\quad (t,x)\in [0,T]\times\R^d,~a,a'\in A\;.
\enqs
\end{Proposition}

\vspace{2mm}

To prove this result we use the following lemma. Observe by definition \reff{relvvartheta} of $v^\pi$ that  
the function $v^\pi$ does not depend on $a$ on the  grid times $\pi$, and we shall denote by misuse of notation: $v^\pi(t_k,x)$, for 
$k$ $\leq$ $n$, $x$ $\in$  $\R^d$. 

\begin{Lemma} \label{lemvvk}
There exists a constant $C$ (not depending on $\pi$) such that
\beqs
|\vartheta^\pi(t,x,a)-v^\pi(t_{k+1},x)| & \leq &  C(1+|x|) |\pi|^{1\over 2}
\enqs
for all  $k=0,\ldots,n-1$, $t\in[t_k,t_{k+1})$, $(x,a)\in\R^d\times A$.
\end{Lemma}
\textbf{Proof.} Fix  $k=0,\ldots,n-1$, $t\in[t_k,t_{k+1})$ and $(x,a)\in\R^d\times A$. Let $(\tilde\Yc,\tilde\Zc,\tilde\Uc)$ be the solution to 
the BSDE
\beqs
\tilde\Yc_s & = & v^\pi(t_{k+1},X^{t,x,a}_{t_{k+1}})+\int_s^{t_{k+1}}f(X_s^{t,x,a},I^{t,a}_s,\tilde\Yc_s,\tilde\Zc_s)ds \\
 & & \; -\int_s^{t_{k+1}} \tilde\Zc_s dW_s  -\int_s^{t_{k+1}}\int_A \tilde\Uc_s(a')\tilde \mu (ds,da')\;,\qquad s\in[t,t_{k+1}]\;.
\enqs
From Proposition \ref{exist uniq EDSR disc const}, Markov property and uniqueness of a solution to the  BSDE \reff{BSDEDC1}-\reff{BSDEDC2}-\reff{BSDEDC3} we have: $\tilde\Yc_s$ $=$ $\vartheta^\pi(s,X_s^{t,x,a},I_s^{t,a})$, for $s$ $\in$ $[t,t_{k+1}]$, and so:
\beq
|\vartheta^\pi(t,x,a)-v^\pi(t_{k+1},x)| &=& \big|\tilde\Yc_t - v^\pi(t_{k+1},x) \big| \nonumber \\
& \leq & \E|v^\pi(t_{k+1},X^{t,x,a}_{t_{k+1}})-v^\pi(t_{k+1},x)|  \nonumber \\
& & + \; \E \Big[ \int_t^{t_{k+1}}\big|f(X^{t,x,a}_s,I^{t,a}_s,\tilde\Yc_s,\tilde\Zc_s)\big|ds \Big].  \label{decompholder v pi}
\enq
From Corollary  \ref{prop stab lip},  we have
\beq \label{decompv cond term holder}
\E|v^\pi(t_{k+1},X^{t,x,a}_{t_{k+1}})-v^\pi(t_{k+1},x)| & \leq & C\sqrt{\E[|X^{t,x,a}_{t_{k+1}}-x|^2]} \; \leq \;  C\sqrt{|\pi|}\;.
\enq
Moreover, by the growth linear condition on $f$ in {\bf (H2)}, and on $\vartheta^\pi$ in Lemma \ref{propcroisbarvn},  we have
\beqs
\E \Big[  \int_t^{t_{k+1}}\big|f(X_s,I_s,\tilde\Yc_s,\tilde\Zc_s)\big|ds \Big]  & \leq &  
C \E \Big[ \int_t^{t_{k+1}}\Big(1+ |X^{t,x,a}_s|  +  |\tilde\Zc_s|\Big) ds \Big].
\enqs
By classical estimates, we have
\beqs
\sup_{s\in[t,T]}\E\big[|X^{t,x,a}_s|^2\big] & \leq & C(1+|x|^2).
\enqs
Moreover, under {\bf (H1)} and {\bf (H2)}, we know from Proposition 4.2 in \cite{BE08} that  there exists a constant $C$ depending only on the Lipschitz constants of $b$, $\sigma$ $f$ and $v^\pi(t_{k+1},.)$ such that
\beqs
\E\big[ \sup_{s\in[t_k,t_{k+1}]} |\tilde\Zc_s|^2 \big] & \leq &   C(1+|x|^2). 
\enqs
This proves that
\beqs
\E \Big[  \int_t^{t_{k+1}}\big|f(X_s,I_s,\tilde\Yc_s,\tilde\Zc_s)\big|ds \Big]  & \leq &  C(1+|x|) |\pi|\;.
\enqs
Combining this last estimate with \reff{decompholder v pi} and \reff{decompv cond term holder}, we get the result
\ep

\vspace{3mm}

\ni \textbf{Proof of Proposition \ref{prop limv indepa}.} The function  $\vartheta$ is l.s.c.  as the supremum of the l.s.c. functions $\vartheta^\pi$. 
Fix $(t,x)\in[0,T)\times \R^d$ and $a,a'\in A$.
Let $(\pi^p)_p$ be a sequence of subdivisions of $[0,T]$ such that $|\pi^p|\downarrow0$ as $p\uparrow\infty$. We define the  sequence  $(t_p)_p$ of $[0,T]$ by
\beqs
t_p & = & \min\big\{ s\in \pi^p~:~s>t\big\}\;,\quad p\geq 0\;.
\enqs
Since $|\pi^p|\rightarrow0$ as $p\rightarrow\infty$ we get $t_p\rightarrow t$ as $p\rightarrow+\infty$. We then have from the previous lemma: 
\beqs
|\vartheta^{\pi^p}(t,x,a)-\vartheta^{\pi^p}(t,x,a')| & \leq & |\vartheta^{\pi^p}(t,x,a)-v^{\pi^p}(t_p,x)|+|v^{\pi^p}(t_p,x) - \vartheta^{\pi^p}(t,x,a')|\\
 & \leq & 2 C |\pi^p|^{1\over 2}\;.
\enqs
Sending $p$ to $\infty$ we obtain that $\vartheta(t,x,a)=\vartheta(t,x,a')$.
\ep

\vspace{3mm}

\begin{Corollary} \label{pointwise}
We  have the identification: $\vartheta$ $=$ $v$, and  the sequence $(v^\pi)_\pi$ also converges to $v$.
\end{Corollary}
\textbf{Proof.}  We proceed in two steps.

\noindent \textbf{Step 1.} \textit{The function $\vartheta$ is a supersolution to \reff{HJB}.} 
\noindent  Since $\vartheta ^{\pi_k}(T,.)=g$ for all $k\geq 1$, we first notice that  $\vartheta (T,.)=g$. 
Next, since $\vartheta$ does not depend on the variable $a$, we have
\beqs
\vartheta ^\pi(t,x,a) & \uparrow & \vartheta(t,x)\quad \mbox{ as }~ |\pi| ~\downarrow~0
\enqs
for any $(t,x,a)\in [0,T]\times\R^d\times A$.  Moreover, since the function $\vartheta$ is l.s.c,  we  have 
\beq
\vartheta \; = \; \vartheta_* & = & \liminf_{|\pi| \rightarrow 0}{}_* \vartheta^\pi,
\label{vinf1}
\enq
 where 
 \beqs
\liminf_{|\pi| \rightarrow 0}  {}_*\vartheta^\pi(t,x,a) \; :=  \liminf_{\tiny{\begin{array}{c}|\pi| \rightarrow 0 \\ (t',x',a')\rightarrow (t,x,a)\\
t'<T\end{array}}} \vartheta^\pi(t',x',a'),\qquad (t,x,a)\in[0,T]\times\R^d\times \R^q\;.
\enqs
Fix now some $(t,x)\in [0,T]\times\R^d$ and $a\in A$ and $(p,q,M)\in \bar J^{2,-}\vartheta (t,x)$, the limiting parabolic subjet of $\vartheta$ at $(t,x)$ (see definition in \cite{craishlio92}). 
From standard stability results, there exists a sequence  $(\pi_k,t_k,x_k,a_k,p_k,q_k,M_k)_k$ such that
\beqs
(p_k,q_k,M_k) & \in & \bar J^{2,-}\vartheta^{\pi_k} (t_k,x_k,a_k)
\enqs
for all $k\geq 1$ and 
\beqs
(t_k,x_k,a_k,\vartheta^{\pi_k} (t_k,x_k,a_k)) & \longrightarrow & (t,x,a,\vartheta (t,x,a))\quad \mbox{ as } k\rightarrow\infty, \; |\pi_k| \rightarrow 0. 
\enqs
From the viscosity supersolution property of $\vartheta^{\pi_k}$ to \reff{seqIPDE} in terms of subjets, we have 
\beqs
-p_k-b(x_k,a_k). q_k -  \frac{1}{2}{\rm tr}(\sigma\sigma\trans(x_k,a_k)M_k)-f\big(x_k,a_k,\vartheta^{\pi_k}(t_k,x_k,a_k),\sigma\trans(x_k,a_k)q_k) & & \\
-\int_A\big(\vartheta^{\pi_k}(t_k,x_k,a')-\vartheta^{\pi_k}(t_k,x_k,a_k)\big)\lambda(da') & \geq & 0
\enqs
for all $k\geq 1$. Sending $k$ to infinity and using \reff{vinf1}, we get
\beqs
-p-b(x,a). q -  \frac{1}{2}{\rm tr}(\sigma\sigma\trans(x,a)M)-f\big(x,a,\vartheta (t,x),\sigma\trans(x,a)q) & \geq & 0.
\enqs
Since $a$ is arbitrary in $A$, this shows 
\beqs
-p- \sup_{a\in A} \big[ b(x,a). q +  \frac{1}{2}{\rm tr}(\sigma\sigma\trans(x,a)M) + f\big(x,a,\vartheta (t,x),\sigma\trans(x,a)q) \big]& \geq & 0,
\enqs
i.e.  the viscosity supersolution property of $\vartheta$ to \reff{HJB}.

\vspace{2mm}

\noindent \textbf{Step 2.} \textit{Comparison.} Since the PDE \reff{HJB} satisfies a comparison principle, we have from the previous step 
$\vartheta \geq v$, and we conclude with  Lemma  \ref{propcroisbarvn} that  $\vartheta= v$. Finally, by definition  \reff{relvvartheta} of $v^\pi$ and from Lemma \ref{propcroisbarvn}, 
we clearly have  $\vartheta^\pi$ $\leq$ $v^\pi$ $\leq$ $v$, which also  proves that $(v^\pi)_\pi$ converges to $v$.
\ep

\vspace{3mm}

In terms of the discretely jump-constrained BSDE,  the convergence result is formulated as follows:

\begin{Proposition} \label{propconv1}
We have $\Yc_t^\pi$ $\leq$ $Y_t^\pi$ $\leq$ $Y_t$, $0\leq t\leq T$, and 
\beqs
\E\Big[\sup_{t\in[0,T]} |Y_t-\Yc_t^\pi|^2\Big]+\E\Big[\sup_{t\in[0,T]}|Y_t-Y_t^\pi|^2\Big]+\E\Big[\int_{0}^T|Z_t-\Zc_t^\pi|^2dt\Big] 
& \rightarrow & 0,
\enqs
as $|\pi|$ goes to zero. 
\end{Proposition}
{\bf Proof.}  Recall from \reff{Yv} and \reff{identif BSDE IPDE} that we have the representation:
\beq \label{idenY}
Y_t \; = \; v(t,X_t), \;\;\; Y_t^\pi \; = \; v^\pi(t,X_t,I_t), \;\; \Yc_t^\pi \; = \; \vartheta(t,X_t,I_t),
\enq
and the first assertion of Lemma \reff{propcroisbarvn}, we clearly have:  $\Yc_t^\pi$ $\leq$ $Y_t^\pi$ $\leq$ $Y_t$, $0\leq t\leq T$. 
The convergence in $\Sc^2$ for $\Yc^\pi$ to $Y$ and $Y^\pi$ to $Y$ comes from the above representation \reff{idenY}, 
the pointwise convergence of  $\vartheta^\pi$ and $v^\pi$ to $v$ in Corollary \ref{pointwise}, and the dominated convergence theorem by recalling that 
$0$ $\leq$ $(v-v^\pi)(t,x,a)$ $\leq$ $(v-\vartheta^\pi)(t,x,a)$ $\leq$ $v(t,x)$ $\leq$ $C(1+|x|)$. 
Let us now turn to the component $Z$.  By definition \reff{BSDEDC1}-\reff{BSDEDC2}-\reff{BSDEDC3} of the discretely jump-constrained BSDE we notice that $\Yc^\pi$ can be written on $[0,T]$ as: 
\beqs
\Yc^\pi_t & = & g(X_T)+\int_t^Tf(X_s,I_s,\Yc_s^\pi,\Zc_s^\pi)-\int_t^T\Zc_s^\pi dW_s-\int_t^T\int_A\Uc_s^\pi(a) \tilde \mu(ds,da)+\Kc^\pi_T-\Kc^\pi_t,
\enqs  
where $\Kc^\pi$ is the nondecreasing process defined by: $\Kc_t^\pi$ $=$ $\sum_{t_k\leq t}(Y_{t_k}^\pi-\Yc_{t_k}^\pi)$, for $t$ $\in$ $[0,T]$. 
Denote by $\delta Y=Y-\Yc^\pi$, $\delta Z=Z-\Zc^\pi$, $\delta U=U-\Uc^\pi$ and $\delta K=K-\Kc^\pi$. 
From It\^o's formula, Young Inequality and  {\bf (H2)}, there exists a constant $C$ such that
\beq \nonumber
& & \E\Big[|\delta Y_t|^2\Big] + {1\over 2}\E\Big[\int_t^T|\delta Z_s|^2ds\Big] +  \frac{1}{2}  \E\Big[\int_t^T|\delta U_s(a)|^2\lambda(da)ds\Big]  \\
& \leq & C\int_t^T\E\Big[|\delta Y_s|^2\Big]ds   + {1\over \eps} \E\Big[\sup_{s\in[0,T]} |\delta Y_s|^2\Big] + \eps\E\Big[\big| \delta K_T-\delta K_t \big|^2\Big]
\label{estim ecar BSDE disc const}
\enq
for all $t\in[0,T]$, with $\eps$ a constant to be chosen later. From the definition of $\delta K$ we have
\beqs
\delta K_T-\delta K_t & = & \delta Y_t - \int_t^T\big(f(X_s,I_s,Y_s,Z_s)-f(X_s,I_s,\Yc^\pi_s,\Zc^\pi_s)\big)ds\\
 & & +\int_0^T\delta Z_s dW_s+\int_t^T\int_A\delta U_s(a)\tilde \mu(ds,da)\;.
\enqs
Therefore, by {(H2)}, we get  the existence of a constant $C'$ such that
\beqs
\E\Big[\big| \delta K_T-\delta K_t\big|^2\Big] & \leq & C'\Big( \E\Big[\sup_{s\in[0,T]} |\delta Y_s|^2\Big] + \E\Big[\int_t^T|\delta Z_s|^2ds\Big]+\E\Big[\int_t^T|\delta U_s(a)|^2\lambda(da)ds\Big] \Big)
\enqs
Taking $\eps={C'\over 4}$  and plugging this last inequality in \reff{estim ecar BSDE disc const}, we get the existence of a constant $C''$ such that 
\beq \label{ZY}
\E\Big[\int_t^T|\delta Z_s|^2ds\Big]+\E\Big[\int_t^T|\delta U_s(a)|^2\lambda(da)ds\Big] & \leq & C''\Big(\E\Big[\sup_{s\in[0,T]} |\delta Y_s|^2\Big]\Big),
\enq
which shows the $L^2(W)$ convergence of $\Zc^\pi$ to $Z$ from the $\Sc^2$ convergence of $\Yc^\pi$ to $Y$. 
\ep

\subsection{Rate of convergence}

We next provide an error estimate for the convergence of the discretely jump-constrained BSDE. We shall combine BSDE methods and  PDE arguments adapted from the shaking coefficients approach of 
Krylov \cite{kr00} and  switching systems approximation of Barles, Jacobsen \cite{barjac07}.  We make further assumptions:

\vspace{2mm}

\ni {\bf (H1')} The functions $b$ and $\sigma$ are uniformly bounded: 
\beqs
\sup_{x\in\R^d,a\in A} |b(x,a)|+|\sigma(x,a)| & < & \infty.
\enqs

\vspace{1mm}

\ni {\bf (H2')} The function $f$ does not depend on $z$: $f(x,a,y,z)$ $=$ $f(x,a,y)$ for all $(x,a,y,z)\in\R^d\times A\times\R\times\R^d$ and
\begin{enumerate}[(i)]
\item the functions $f(.,.,0)$ and $g$ are uniformly bounded: 
\beqs
\sup_{x\in\R^d,a\in A} |f(x,a,0)| + |g(x)| & < & \infty,
\enqs
\item for all $(x,a)$ $\in$ $\R^d\times A$, $y$ $\mapsto$ $f(x,a,y)$ is convex.
\end{enumerate}

\vspace{2mm}

Under these assumptions,  we obtain the rate of convergence for $v^\pi$ and $\vartheta^\pi$ towards  $v$.

\begin{Theorem}\label{thm convdisc ctr}
Under {\bf (H1')} and {\bf (H2')},  there exists a constant $C$ such that
\beqs
0~~\leq~~ v(t,x)-v ^\pi(t,x,a) ~~\leq~~ v(t,x)-\vartheta ^\pi(t,x,a) & \leq & C |\pi|^{1\over 10}
\enqs
for all $(t,x,a)\in [0,T]\times\R^d\times A$ and all grid $\pi$ with $|\pi|$ $\leq$ $1$.  
Moreover, when $f(x,a)$ does not depend on $y$,  the rate of convergence is improved to $|\pi|^{1\over 6}$.
\end{Theorem}

Before proving this result, we give  as corollary the  rate of convergence  for the discretely jump-constrained BSDE.  
 
\begin{Corollary}\label{cormainTH}
Under {\bf (H1')} and {\bf (H2')},  there exists a constant $C$ such that
\beqs
\E\Big[\sup_{t\in[0,T]}|Y_t-\Yc_t^\pi|^2\Big]+\E\Big[\sup_{t\in[0,T]}|Y_t-Y_t^\pi|^2\Big]+\E\Big[\int_{0}^T|Z_t-\Zc_t^\pi|^2dt\Big] & \leq & C |\pi|^{1\over 5}. 
\enqs
for all grid $\pi$ with $|\pi|$ $\leq$ $1$, and the above rate is improved to $|\pi|^{1\over 3}$ when $f(x,a)$ does not depend on $y$. 
\end{Corollary}
\textbf{Proof.} From  the representation \reff{idenY},  and Theorem \ref{thm convdisc ctr},  we immediately have
\beq\label{estimdelta Y}
\E\Big[\sup_{t\in[0,T]}|Y_t-\Yc_t^\pi|^2\Big]+\E\Big[\sup_{t\in[0,T]}|Y_t-Y_t^\pi|^2\Big] & \leq & C |\pi|^{1\over 5},
\enq
(resp. $|\pi|^{1\over 3}$ when $f(x,a)$ does not depend on $y$). 
Finally, the convergence rate for $Z$ follows from the inequality \reff{ZY}. 
\ep

 \begin{Remark}
 {\rm   The above  convergence rate $|\pi|^{1\over 10}$   is the  optimal rate that  one can prove  in our generalized stochastic control context with fully  nonlinear HJB equation by PDE approach and shaking coefficients  technique, see \cite{kr00}, \cite{barjac07}, \cite{ftw11} or \cite{tan13}.  However,  this rate  may not be the sharpest one.  In the case of  continuously reflected BSDEs, i.e. BSDEs with upper or lower constraint on $Y$,    it is known  that $Y$ can be approximated  by discretely reflected BSDEs, i.e.  BSDEs where reflection on $Y$ operates a finite set of times  on the grid $\pi$, with a rate $|\pi|^{1\over 2}$ 
 (see \cite{balpag03}). The standard arguments for proving this rate is based on the representation of the continuously (resp. discretely) reflected BSDE as optimal stopping problems where stopping is  possible over the whole interval time (resp.  only on the grid times).  In our jump-constrained case,  we know from \cite{kp12} that  the minimal solution to the BSDE with nonpositive jumps has the  stochastic control  representation  \reff{constov} when $f(x,a)$ does not depend on $y$ and $z$. We shall  prove  an  analog representation  for discretely jump-constrained BSDEs, and this helps to improve the rate of convergence from $|\pi|^{1\over 10}$ to $|\pi|^{1\over 6}$. 
 \ep
 }
 \end{Remark}

\vspace{2mm}

The rest of this section is devoted to the proof of Theorem \ref{thm convdisc ctr}.  We first consider the special case where $f(x,a)$ does not depend on $y$, and then address the case $f(x,a,y)$. 

\vspace{1mm}

\noindent {\bf Proof of Theorem \ref{thm convdisc ctr} in the case $f(x,a)$}. 

\ni In the case where $f(x,a)$ does not depend on $y,z$, by  (linear) Feynman-Kac formula for $\vartheta^\pi$ solution to \reff{seqIPDE}, and by definition of $v^\pi$ in \reff{relvvartheta}, we have: 
\beqs
v^\pi(t_k,x) &=& \sup_{a \in A} \E\Big[ \int_{t_k}^{t_{k+1}}  f(X_t^{t_k,x,a},I_t^{t_k,a}) dt + v^\pi(t_{k+1},X_{t_{k+1}}^{t_k,x,a}) \Big], \;\; k  \leq n-1, \; x \in \R^d. 
\enqs
By induction, this  dynamic programming relation leads to the following stochastic control problem with  discrete time policies:  
\beqs
v^\pi(t_k,x) &=& \sup_{\alpha\in\Ac_{\F}^\pi} \E\Big[ \int_{t_k}^T f(\bar X_t^{t_k,x,\alpha},\bar I_t^\alpha) dt + g(\bar X_T^{t_k,x,\alpha}) \Big],
\enqs
where $\Ac_\F^\pi$ is the set of  discrete time processes $\alpha$ $=$ $(\alpha_{t_j})_{j\leq n-1}$, with $\alpha_{t_j}$ $\Fc_{t_j}$-measurable, valued in $A$, and 
\beqs
\bar X_t^{t_k,x,\alpha} &=& x  +  \int_{t_k}^t b(\bar X_s^{t,x,\alpha},\bar I_s^\alpha) ds + \int_{t_k}^t \sigma(\bar X_s^{t_k,x,\alpha},\bar I_s^\alpha) dW_s, \;\;\; t_k \leq t \leq T, \\
\bar I_t^\alpha &=&   \alpha_{t_j} + \int_{(t_j,t]} \int_A (a- \bar I_{s^-}^\alpha) \mu(ds,da), \;\;\; t_j \leq t < t_{j+1},  \; j \leq n-1.
\enqs
In other words, $v^\pi(t_k,x)$ corresponds to the value function  for a stochastic control problem where the controller can act only at the dates $t_j$ of the grid $\pi$, and then let the regime of the coefficients of the diffusion evolve  according to the Poisson random measure $\mu$.  Let us introduce the following stochastic control problem with piece-wise constant control policies: 
\beqs
\tilde v^\pi(t_k,x) &=& \sup_{\alpha\in\Ac_{\F}^\pi} \E\Big[ \int_{t_k}^T f(\tilde X_t^{t_k,x,\alpha},\tilde I_t^\alpha) dt + g(\tilde X_T^{t_k,x,\alpha}) \Big],
\enqs
where for $\alpha$ $=$ $(\alpha_{t_j})_{j\leq n-1}$ $\in$ $\Ac_\F^\pi$: 
\beqs
\tilde X_t^{t_k,x,\alpha} &=& x  +  \int_{t_k}^t b(\tilde X_s^{t,x,\alpha},\tilde I_s^\alpha) ds + \int_{t_k}^t \sigma(\tilde X_s^{t_k,x,\alpha},\tilde I_s^\alpha) dW_s, \;\;\; t_k \leq t \leq T, \\
\tilde I_t^\alpha &=&   \alpha_{t_j}, \;\;\;\;\;\;\;   t_j \leq t < t_{j+1},  \; j \leq n-1.
\enqs
It is shown in \cite{kr99} that $\tilde v^\pi$ approximates the value function $v$  for the controlled diffusion problem \reff{constov},  solution to the HJB equation \reff{HJB}, with a rate $|\pi|^{1\over 6}$:
\beq \label{approkrylov}
0 \; \leq \; v(t_k,x) - \tilde v^\pi(t_k,x) & \leq & C |\pi|^{1\over 6},
\enq
for all $t_k$ $\in$ $\pi$, $x$ $\in$ $\R^d$.  Now, recalling that $A$ is compact and $\lambda(A)$ $<$ $\infty$,  it is clear that there exists some positive constant $C$ such that for all 
$\alpha$ $\in$ $\Ac_\F^\pi$, $j$ $\leq$ $n-1$: 
\beqs
\E \Big[ \sup_{t \in [t_j,t_{j+1})} |\bar I_t^\alpha -  \tilde I_t^\alpha|^2 \Big] & \leq & C |\pi|,   
\enqs
and then by standard arguments under Lipschitz condition on $b$, $\sigma$: 
\beqs
\E \Big[ \sup_{ t\in [t_j,t_{j+1}]} |\bar X_t^{t_k,x,\alpha}-\tilde X_t^{t_k,x,\alpha}|^2 \Big] & \leq & C|\pi|, \;\;\; k \leq j \leq n-1, \; x \in \R^d.     
\enqs
By  the Lipschitz  conditions on $f$ and $g$, it follows that
\beqs
|v^\pi(t_k,x) - \tilde v^\pi(t_k,x) | & \leq & C |\pi|^{1\over 2}, 
\enqs
and thus with \reff{approkrylov}: 
\beqs
0 \; \leq \; \sup_{x\in\R^d} (v-v^\pi)(t_k,x) & \leq & C |\pi|^{1\over 6}. 
\enqs
Finally, by combining with the estimate in Lemma \ref{lemvvk}, which gives  actually under {\bf (H2')}(i): 
\beqs
|\vartheta^\pi(t,x,a) - v^\pi(t_{k+1},x)| & \leq & C|\pi|^{1\over 2}, \;\;\; t \in [t_k,t_{k+1}), (x,a) \in \R^d\times A,
\enqs
together with the $1/2$-H\"older property of $v$ in time (see \reff{vlip}),  we obtain:
\beqs
\sup_{(t,x,a)\in[0,T]\times\R^d\times A} (v - \vartheta^\pi)(t,x,a) & \leq & C (|\pi|^{1\over 6} + |\pi|^{1\over 2}) \; \leq \; C |\pi|^{1\over 6},
\enqs
for $|\pi|$ $\leq$ $1$. This ends the proof. 
\ep

\vspace{3mm}

Let us now turn to the case where $f(x,a,y)$ may also depend on $y$. We cannot rely anymore on the convergence  rate  result in \cite{kr99}. 
Instead, recalling that $A$ is compact and  since $\sigma$, $b$ and $f$ are Lipschitz in $(x,a)$, we shall  apply  the switching 
system method of Barles and Jacobsen \cite{barjac07}, which is a variation of the shaken coefficients method and smoothing technique used in Krylov \cite{kr00}, in order to obtain approximate smooth subsolution to 
\reff{HJB}.  By Lemmas  3.3 and 3.4 in \cite{barjac07},  one can find a family of  smooth functions $(w_\eps)_{0<\eps\leq 1}$ on $[0,T]\times\R^d$ such that: 
\beq
\sup_{[0,T]\times\R^d}  |w_\eps| & \leq & C,   \label{wepsbor} \\
\sup_{[0,T]\times\R^d} |w_\eps-w| & \leq & C \eps^{1\over 3},   \label{wepsapprox} \\
 \sup_{[0,T]\times\R^d} |\partial_t^{\beta_0}D^\beta  w_\eps| & \leq &  C \eps^{1-2\beta_0-\sum_{i=1}^d \beta^i},  \;\;  \beta_0\in \N, \; \beta=(\beta^1,\ldots,\beta^d)\in \N^d, \label{cond-der-w-eps}
 \enq 
for some positive constant $C$  independent of $\eps$, and by convexity of $f$ in {\bf (H2')}(ii),  for any $\eps$ $\in$ $(0,1]$, $(t,x)$ $\in$ $[0,T]\times\R^d$, there exists  $a_{t,x,\eps}$ $\in$ $A$ such that: 
\beq \label{sousolweps}
- \Lc^{a_{t,x,\eps}} w_\eps(t,x)  - f(x,a_{t,x,\eps},w_\eps(t,x)) & \geq & 0. 
\enq

\vspace{2mm}

Recalling the  definition of the operator $\T_\pi^k$ in \reff{operaT}, we define for any function $\varphi$ on $[0,T]\times\R^d\times A$, Lipschitz in $(x,a)$:  
\beqs
\T_\pi[\varphi](t,x,a) &:=& \T_\pi^k[\varphi(t_{k+1},.,.)](t,x,a), \;\;\;\; t \in [t_k,t_{k+1}),  \;  (x,a) \in \R^d\times A, 
\enqs
for $k$ $=$ $0,\ldots,n-1$, and  
\beqs
\mathbb{S}_\pi[\varphi](t,x,a) & := & {1\over |\pi|}\Big[\varphi(t,x)-\T_\pi[\varphi](t,x,a)  \\
& & \;\;\; \;\;\;\;  +  (t_{k+1}-t) \big( \Lc^a \varphi(t,x)+f(x,a,\varphi(t,x) \big) \Big],
\enqs
for $(t,x,a)$ $\in$ $[t_k,t_{k+1})\times\R^d\times A$, $k$ $\leq$ $n-1$.

\vspace{1mm}

We have the following key error bound on $\S_\pi$.

\begin{Lemma}\label{lem-Sn}
Let {\bf (H1')} and {\bf (H2')}(i) hold.  There exists a constant $C$  such that
\beqs
|\mathbb{S}_\pi[\varphi_\eps](t,x,a)| & \leq & C \Big(|\pi|^{1\over 2}(1+\eps^{-1})+ |\pi|{\eps^{-3}} \Big), \;\;\; (t,x,a) \in [0,T]\times\R^d\times A,
\enqs
for any family $(\varphi_\eps)_\eps$ of smooth functions on $[0,T]\times\R^d$ satisfying  \reff{wepsbor} and \reff{cond-der-w-eps}. 
\end{Lemma}
\textbf{Proof.} Fix  $(t,x,a)\in[0,T]\times\R^d\times A$. If $t=T$, we have $|\mathbb{S}_\pi[\varphi_\eps](t,x,a)|=0$. Suppose that $t <T$ and fix $k\leq n$ such that  $t\in [t_k,t_{k+1})$. 
Given a smooth  function $\varphi_\eps$ satisfying \reff{wepsbor} and \reff{cond-der-w-eps}, we split: 
\beqs
|\mathbb{S}_\pi[\varphi_\eps](t,x,a)| & \leq & A_\eps(t,x,a)+ B_\eps(t,x,a),
\enqs
where
\beqs
A_\eps(t,x,a) & := & {1\over |\pi|}\Big|\T_\pi[\varphi_\eps](t,x,a)-\E\big[ \varphi_\eps({t_{k+1}},X_{t_{k+1}}^{t,x,a})\big]- (t_{k+1}-t)f(x,a,\varphi_\eps(t,x)\big)\Big|\;,
\enqs
and
\beqs
B_\eps(t,x,a) & := & {1\over |\pi|}\Big|\E\big[ \varphi_\eps({t_{k+1}},X_{t_{k+1}}^{t,x,a}) \big] - \varphi_\eps(t,x)-(t_{k+1}-t)
\Lc^a\varphi_\eps(t,x)  \Big|,
\enqs
and we study  each term $A_\eps$ and $B_\eps$ separately. 

\vspace{2mm}

\ni \textbf{1.} \textit{Estimate on $A_\eps(t,x,a)$.}

\ni Define $(Y^{\varphi_\eps},Z^{\varphi_\eps},U^{\varphi_\eps})$ as the solution to  the BSDE on $[t,t_{k+1}]$: 
\beq
Y_s^{\varphi_\eps}  & = & \varphi_\eps({t_{k+1}},X_{t_{k+1}}^{t,x,a}) + 
\int_s^{t_{k+1}}  f(X_{r}^{t,x,a},I_{r}^{t,a},Y_r^{\varphi_\eps})  dr  \nonumber \\
 & & \; - \int_s^{t_{k+1}} Z_r^{\varphi_\eps} dW_r - \int_s^{t_{k+1}}\int_A U_r^{\varphi_\eps}(a)\tilde \mu(dr,da)\;,\quad s \in [t,t_{k+1}].  \label{Yvarphi}
\enq
From Theorems  3.4 and 3.5 in \cite{bbp97}, we have $Y_t^{\varphi_\eps}$ $=$ $\T_\pi[\varphi_\eps](t,x,a)$, and by taking expectation in 
\reff{Yvarphi}, we thus get:
\beqs
Y_t^{\varphi_\eps} \; = \;  \T_\pi[\varphi_\eps](t,x,a) &=& \E\big[  \varphi_\eps({t_{k+1}},X_{t_{k+1}}^{t,x,a})  \big]  + 
\E \Big[ \int_t^{t_{k+1}}  f(X_{s}^{t,x,a},I_{s}^{t,a},Y_s^{\varphi_\eps})  ds  \Big]
\enqs
and so: 
\beqs
A_\eps(t,x,a) & \leq &  {1\over |\pi|}  \E \Big[  \int_t^{t_{k+1}}  \big|  f(X_{s}^{t,x,a},I_{s}^{t,a},Y_s^{\varphi_\eps}) - f(x,a,\varphi_\eps(t,x))\big| ds \Big] \\
& \leq & C \Big( \E\big[ \sup_{s\in [t,t_{k+1}]} |X_s^{t,x,a}-x| + |I_s^{t,a}-a| \big]   +  \E\big[\sup_{s\in [t,t_{k+1}]} |Y_s^{\varphi_\eps}-\varphi_\eps(t,x)| \big] \Big),
\enqs
by the Lipschitz continuity of $f$.  From standard estimate for SDE, we have (recall that the coefficients $b$ and $\sigma$ are bounded under {\bf (H1')} and $A$ is compact):
\beq \label{estimXI}
\E\big[ \sup_{s\in [t,t_{k+1}]} |X_s^{t,x,a}-x| + |I_s^{t,a}-a| \big]  & \leq & C |\pi|^{1\over 2}. 
\enq
Moreover, by \reff{Yvarphi}, the boundedness condition in {\bf (H2')}(i) together with the Lipschitz condition of $f$, and Burkholder-Davis-Gundy inequality, we have:
\beqs
\E \Big[ \sup_{s\in [t,t_{k+1}]} |Y_s^{\varphi_\eps}-\varphi_\eps(t,x)|  \Big] & \leq &
 \E \big[ | \varphi_\eps({t_{k+1}},X_{t_{k+1}}^{t,x,a}) - \varphi_\eps(t,x) | \big]  \\
 & & \;\; + \;  C |\pi| \E\big[ 1 + \sup_{s \in[t,t_{k+1}]}  |Y_s^{\varphi_\eps}| \big] \\
 & & \;\ + \; C|\pi| \Big( \E \big[ \sup_{s \in[t,t_{k+1}]}|Z_s^{\varphi_\eps}|^2] +  \E \big[ \sup_{s \in[t,t_{k+1}]} \int_A |U_s^{\varphi_\eps}(a)|^2 \lambda(da)] \Big). 
\enqs
From standard estimate for the BSDE \reff{Yvarphi}, we have:
\beqs
\E\big[\sup_{s\in [t,t_{k+1}]} |Y_s^{\varphi_\eps}|^2 \big] & \leq & C,
\enqs
for some positive constant $C$ depending only on the Lipschitz constant of $f$, the upper bound of $|f(x,a,0,0)|$ in {\bf (H2')}(i), and the upper bound of $|\varphi_\eps|$ in \reff{wepsbor}.  
Moreover,  from the estimate in Proposition 4.2 in \cite{BE08} about the coefficients $Z^{\varphi_\eps}$ and $U^{\varphi_\eps}$ of the BSDE  with jumps \reff{Yvarphi},  there exists some constant $C$ depending only on the Lipschitz constant of $b,\sigma,f$, and of the Lipschitz constant  of $\varphi_\eps(t_{k+1},.)$ (which does not depend on $\eps$ by \reff{cond-der-w-eps}), such that:
\beqs
\E \big[ \sup_{s \in[t,t_{k+1}]}|Z_s^{\varphi_\eps}|^2] +  \E \big[ \sup_{s \in[t,t_{k+1}]} \int_A |U_s^{\varphi_\eps}(a)|^2 \lambda(da)] 
& \leq & C. 
\enqs
From \reff{cond-der-w-eps},  we then have:
\beqs
\E \Big[ \sup_{s\in [t,t_{k+1}]} |Y_s^{\varphi_\eps}-\varphi_\eps(t,x)|  \Big] & \leq & C \big( |t_{k+1}-t|\eps^{-1}
+ \E\big[ | X_{t_{k+1}}^{t,x,a} - x|\big]  + |\pi|  \big) \\
& \leq &  C |\pi|^{1\over 2} \big(1 + \eps^{-1} \big),
\enqs
by \reff{estimXI}.   This leads to the error bound for $A_\eps(t,x,a)$:
\beqs
A_\eps(t,x,a) & \leq &  C |\pi|^{1\over 2} (1 + \eps^{-1}). 
\enqs
 
\vspace{2mm}

\ni \textbf{2.} \textit{Estimate on $B_\eps(t,x,a)$.}

\ni From It\^o's formula we have
\beqs
B_\eps(t,x,a) & =  & {1\over |\pi|} \Big| \E\Big[ \int_t^{t_{k+1}} \big(
 \Lc^{I_s^{t,a}}  \varphi_\eps(s,X^{t,x,a}_s)-   \Lc^{a}\varphi_\eps(t,x) \big) ds \Big] \Big| \\
 & \leq &  B_\eps^1(t,x,a) + B_\eps^2(t,x,a)
 \enqs
 where
 \beqs
B_\eps^1(t,x,a) &  = &  {1\over |\pi|}\E\Big[ \int_t^{t_{k+1}}\Big|\big(b(X_s^{t,x,a},I^{t,a}_s\big)-b(x,a)).D_x\varphi_\eps(s,X^{t,x,a}_s)  \\
 & & \hspace{1.3cm} + \;  {1\over 2} {\rm tr}  \big[
 \big(\sigma\sigma\trans(X_s^{t,x,a},I^{t,a}_s)  -\sigma\sigma\trans(x,a)\big)D_x^2\varphi_\eps(t,x)\big]  \Big| ds \Big]
\enqs
and
\beqs
B_\eps^2(t,x,a) & = &  {1\over |\pi|}\E \Big[ \int_t^{t_{k+1}}\big|\tilde \Lc_{t,x}^a\varphi_\eps(s,X_s^{t,x,a}) -\tilde \Lc_{t,x}^a\varphi_\eps(t,x) \big|ds \Big]
 \;,
\enqs
with $\tilde \Lc_{t,x}^a$ defined by
\beqs
\tilde \Lc_{t,x}^a\varphi_\eps(t',x') & = & \Dt{\varphi_\eps}(t',x') + b(x,a). D_x \varphi_\eps(t',x') 
+ \frac{1}{2}{\rm tr}\big(\sigma\sigma\trans(x,a)D_x^2\varphi_\eps(t',x')\big).
\enqs
Under {\bf (H1)}, {\bf (H1')},  and  by \reff{cond-der-w-eps}, we have 
\beqs
B_\eps^1(t,x,a) 
& \leq & C(1+\eps^{-1}) \E\Big[ \sup_{s\in [t,t_{k+1}]} |X^{t,x,a}_s-x| +|I^{t,a}_s-a| \Big] \\
& \leq & C(1+\eps^{-1})|\pi|^{1\over 2},
\enqs
where we used again \reff{estimXI}.  On the other hand, since $\varphi_\eps$ is smooth, we have from It\^o's formula
\beqs
B_\eps^2(t,x,a) & = &  {1\over |\pi|}\E\Big[\int_t^{t_{k+1}}\Big|\int_t^s\Lc^{I_r^{t,a}}\tilde \Lc_{t,x}^a\phi(r,X_r^{t,x,a})dr \Big|ds \Big]\;.
\enqs
Under {\bf (H1')} and  by \reff{cond-der-w-eps}, we then see that
\beqs
B_\eps^2(t,x,a) & \leq & C|\pi|\eps^{-3},
\enqs
and so:
\beqs
B_\eps(t,x,a) & \leq &  C \Big(|\pi|^{1\over 2}(1+\eps^{-1})+ |\pi|{\eps^{-3}} \Big).
\enqs
Together with the estimate for $A_\eps(t,x,a)$, this proves the  error bound for $|\mathbb{S}_\pi[\varphi_\eps](t,x,a)|$. 
\ep

\vspace{3mm}

We next state a maximum principle type result  for  the operator $\T_\pi$.

\begin{Lemma}\label{thm-comp-ipde}
Let $\varphi$ and $\psi$ be two  functions on $[0,T]\times\R^d\times A$, Lipschitz  in $(x,a)$. Then, 
there exists some positive constant $C$ independent of $\pi$ such that
\beqs \label{maxprinTpi}
\sup_{(t,x,a)\in [t_k,t_{k+1}]\times\R^d\times A}(\T_\pi[\varphi]-\T_\pi[\psi])(t,x,a) & \leq & e^{C|\pi|} \sup_{(x,a)\in \R^d\times A}(\varphi-\psi)(t_{k+1},x,a)\;,
\enqs
 for all $k$ $=$ $0,\ldots,n-1$. 
\end{Lemma}
\textbf{Proof.} 
 Fix $k$ $\leq$ $n-1$, and set
\beqs
M &:=& \sup_{(x,a)\in \R^d\times A}(\varphi-\psi)(t_{k+1},x,a).
\enqs
We can assume w.l.o.g. that $M$ $<$ $\infty$ since otherwise the required inequality is trivial. Let us denote by $\Delta v$ the function 
\beqs
\Delta v(t,x,a)  &=&  \T_\pi[\varphi] (t,x,a)-  \T_\pi[\psi] (t,x,a), 
\enqs   
for all $(t,x,a)\in [t_{k},t_{k+1}]\times\R^d\times A$. By definition of $\T_\pi$, and from  the Lipschitz condition of $f$, we see that  $\Delta v$ is a viscosity subsolution to 
\begin{equation}\left\{
\begin{array}{l} -  \Lc^a  \Delta v (t,x,a)- C\big(|\Delta v(t,x,a)|+|D \Delta v(t,x,a)|\big)\\
-\int_A\big(\Delta v(t,x,a')-\Delta v(t,x,a)\big)\lambda(da') = 0, \quad \mbox{ for } \;(t,x,a)\in [t_{k},t_{k+1})\times\R^d\times A, \\
\Delta v (t_{k+1},x,a)\leq M\;,\qquad \mbox{ for } \;(x,a)\in \R^d\times A\;. 
\end{array}\right.
\end{equation}
Then, we easily check that the function $\Phi$ defined by 
\beqs
\Phi(t,x,a) & = & Me^{C(t_{k+1}-t)}\;,\quad (t,x,a)\in [t_k,t_{k+1}]\times \R^d \times A\;,
\enqs
is a solution to 
\begin{equation}\left\{
\begin{array}{l} -  \Lc^a  \Phi (t,x,a)- C\big(|\Phi(t,x,a)|+|D \Phi(t,x,a)|\big)\\
-\int_A\big(\Phi(t,x,a')-\Phi(t,x,a)\big)\lambda(da')  \; =  \; 0, \quad \mbox{ for } \;(t,x,a)\in [t_{k},t_{k+1})\times\R^d\times A, \\
\Phi (t_{k+1},x,a) \; = \;  M\;,\qquad \mbox{ for } \;(x,a)\in \R^d\times A\;. 
\end{array}\right.
\end{equation}
From the comparison theorem in \cite{bbp97} for  viscosity solutions of semi-linear IPDEs, we get that $\Delta v$ $\leq$ $\Phi$ on $[t_k,t_{k+1}]\times\R^d\times A$, which proves the required inequality. 
\ep

\vspace{3mm}

\ni {\bf Proof of Theorem \ref{thm convdisc ctr}.}  By  \reff{relvvartheta} and \reff{varthetapi}, we observe that $v^\pi$ is a fixed point of $\T_\pi$, i.e.
\beqs
\T_\pi[v^\pi] & =  & v^\pi. 
\enqs
On the other hand, by \reff{sousolweps}, and the estimate of Lemma \ref{lem-Sn} applied to $w_\eps$, we have: 
\beqs
w_\eps(t,x) - \T_\pi[w_\eps](t,x,a_{t,x,\eps}) & \leq &  |\pi|\S_\pi[w_\eps](t,x,a_{t,x,\eps})  \; \leq \; C |\pi| \bar S(\pi,\eps) 
\enqs
where we set: $\bar S(\pi,\eps)$ $=$ $(|\pi|^{3\over2}(1+\eps^{-1})+ |\pi|^2{\eps^{-3}})$.  
Fix $k$ $\leq$ $n-1$. By  Lemma \ref{thm-comp-ipde}, we then have for all $t$ $\in$ $[t_k,t_{k+1}]$, $x$ $\in$ $\R^d$:
\beq
w_\eps(t,x)  - v^\pi(t,x,a_{t,x,\eps}) & = & w_\eps(t,x) - \T_\pi[w_\eps](t,x,a_{t,x,\eps}) +  (\T_\pi[w_\eps] - \T_\pi[v^\pi])(t,x,a_{t,x,\eps}) \nonumber  \\
& \leq & C |\pi| \bar S(\pi,\eps)   +    e^{C |\pi|} \sup_{(x,a)\in\R^d\times A}(w_\eps - v^\pi)(t_{k+1},x,a).  \label{interweps}
\enq
Recalling by its very definition that $v^\pi$ does not depend on $a$ $\in$ $A$ on the grid times of $\pi$, and denoting then 
$M_k$ $:=$ $\sup_{x\in\R^d} (w_\eps - v^\pi)(t_{k},x)$,  we have by \reff{interweps} the relation: 
\beqs
M_k & \leq & C |\pi| \bar S(\pi,\eps)   + e^{C|\pi|} M_{k+1}.
\enqs
By induction, this yields:
\beqs
\sup_{x\in\R^d} (w_\eps - v^\pi)(t_{k},x) & \leq & C \frac{e^{Cn|\pi|}-1}{e^{C|\pi|}-1} |\pi| \bar S(\pi,\eps)  +  e^{C n|\pi|} \sup_{x\in\R^d} (w_\eps - v^\pi)(T,x) \\
& \leq &  C\bar S(\pi,\eps) +  C  \sup_{x\in\R^d} (w_\eps - v)(T,x), 
\enqs
since $n|\pi|$ is bounded and $v(T,x)$ $=$ $v^\pi(T,x)$ ($=$ $g(x)$). 
From \reff{wepsapprox}, we then get:
\beqs
\sup_{x\in\R^d} (v - v^\pi)(t_{k},x)  & \leq &  C \big(  \eps^{1\over 3} + |\pi|^{1\over 2}(1+\eps^{-1})+ |\pi|{\eps^{-3}} \big). 
\enqs
By minimizing the r.h.s of this estimate with respect to $\eps$, this leads to the error bound when  taking $\eps$ $=$ $|\pi|^{3\over10}$ $\leq$ $1$: 
\beqs
\sup_{x\in\R^d} (v - v^\pi)(t_{k},x)  & \leq &  C  |\pi|^{1\over 10}.  
\enqs
Finally, by combining with the estimate in Lemma \ref{lemvvk}, which gives  actually under {\bf (H2')}(i): 
\beqs
|\vartheta^\pi(t,x,a) - v^\pi(t_{k+1},x)| & \leq & C|\pi|^{1\over 2}, \;\;\; t \in [t_k,t_{k+1}), (x,a) \in \R^d\times A,
\enqs
together with the $1/2$-H\"older property of $v$ in time (see \reff{vlip}), we obtain:
\beqs
\sup_{(t,x,a)\in[0,T]\times\R^d\times A} (v - \vartheta^\pi)(t,x,a) & \leq & C (|\pi|^{1\over 10} + |\pi|^{1\over 2}) \; \leq \; C |\pi|^{1\over 10}. 
\enqs
This ends the proof.
\ep

\section{Approximation scheme for  jump-constrained BSDE and stochastic control problem}

\setcounter{equation}{0} \setcounter{Assumption}{0}
\setcounter{Theorem}{0} \setcounter{Proposition}{0}
\setcounter{Corollary}{0} \setcounter{Lemma}{0}
\setcounter{Definition}{0} \setcounter{Remark}{0}


We consider the discrete time approximation of  the discretely jump-constrained BSDE  in the case  where $f(x,a,y)$ does not depend on $z$, and define 
the scheme $(\bar Y^\pi,\bar\Yc^\pi,\bar\Zc^\pi)$ by  induction on the grid $\pi$ $=$  $\{t_0=0 < \ldots <t_k < \ldots <t_n = T\}$  by:
\begin{equation}\label{schemeBSDE2}
\left\{ \begin{array}{rcl}
\bar Y_T^\pi \; = \; \bar \Yc_T^\pi &=& g(\bar X_T^\pi) \\
\bar\Yc_{t_k}^\pi &=& \E_{t_k}\big[ \bar Y_{t_{k+1}}^\pi  \big] +   f(\bar X_{t_k}^\pi,I_{t_k},\bar\Yc_{t_k}^\pi )  \Delta t_k  \\
\bar Y_{t_k}^\pi &=& \esssup_{a \in A} \E_{t_k,a} \big[ \bar\Yc_{t_k}^\pi \big], \;\;\; k = 0, \ldots, n-1, 
\end{array}
\right.
\end{equation}
where $\Delta t_k$ $=$ $t_{k+1}-t_k$,  $\E_{t_k}[.]$ stands for $\E[.|\Fc_{t_k}]$, and $\E_{t_k,a}[.]$ for  $\E[.|\Fc_{t_k},I_{t_k}=a]$.  

\vspace{2mm}

By induction argument, we easily see that $\bar\Yc_{t_k}^\pi$ is a deterministic function of  $(\bar X_{t_k}^\pi,I_{t_k})$, while $\bar Y_{t_k}^\pi$ is a deterministic function of $\bar X_{t_k}^\pi$, for $k$ $=$ $0,\ldots,n$, and by the Markov pro\-perty of  the process $(\bar X^\pi,I)$, the conditional expectations in \reff{schemeBSDE2} can be replaced by the corresponding regressions:
\beqs
 \E_{t_k}\big[ \bar Y_{t_{k+1}}^\pi  \big]  \; = \;  \E\big[  \bar Y_{t_{k+1}}^\pi \big| \bar X_{t_k}^\pi, I_{t_k} \big] & \mbox{ and } & 
 \E_{t_k,a}\big[ \bar\Yc_{t_k}^\pi] \; = \; \E\big[ \bar\Yc_{t_k}^\pi \big|  \bar X_{t_k}^\pi, I_{t_k} = a \big].
\enqs
We then have:
\beqs
\bar\Yc_{t_k}^\pi \; = \; \bar\vartheta^\pi_k(\bar X_{t_k}^\pi,I_{t_k}), & & Y_{t_k}^\pi \; = \; \bar v_k^\pi(\bar X_{t_k}^\pi), 
\enqs 
for some sequence of functions $(\bar\vartheta_k^\pi)_k$ and $(\bar v_k^\pi)_k$  defined respectively on $\R^d\times A$ and $\R^d$ by backward induction:
\begin{equation}\label{schemeBSDE3}
\left\{ \begin{array}{rcl}
\bar v_n^\pi(x,a)  \; = \; \bar\vartheta_n^\pi(x) &=& g(x) \\
\bar\vartheta_k^\pi(x,a) &=& \E \Big[ \bar v_{k+1}^\pi(\bar X_{t_{k+1}}^\pi,I_{t_{k+1}})  \big| (\bar X_{t_k}^\pi, I_{t_k}) = (x,a) \Big] +  f(x,a,\bar\vartheta_k^\pi(x,a))  \Delta t_k \\
\bar v_k^\pi (x)  &=& \sup_{a \in A} \bar\vartheta_k^\pi(x,a), \;\;\; k  = 0,\ldots,n-1. 
\end{array}
\right.
\end{equation}
There are well-known different methods (Longstaff-Schwartz least square regression, quantization, Malliavin integration by parts, see e.g. \cite{balpag03}, \cite{gobetal06}, \cite{boutou04}) 
for computing the above conditional expectations, and so the functions $\bar\vartheta_k^\pi$ and 
$\bar v_k^\pi$. It appears that in our context, the simulation-regression method on basis functions defined on $\R^d\times A$, is quite suitable since it allows us to derive at each time step 
$k$ $\leq$ $n-1$, a functional form  $\hat a_k(x)$, which attains the supremum over $A$ in $\bar\vartheta_k^\pi(x,a)$.  We shall see later in this section that the feedback control 
$(\hat a_k(x))_k$ provides an approximation of the optimal control   for the HJB equation associated to a stochastic control  problem when $f(x,a)$ does not depend on $y$. 
We refer to our companion paper  \cite{khalanpha13b}  for the details about the computation of functions $\bar\vartheta_k^\pi$, $\bar v_k^\pi$, $\hat a_k$ by simulation-regression methods, 
and the associated error analysis. 

\subsection{Error estimate for the discrete time scheme}

The main result of this section is to state an error bound between the component  $Y^\pi$ of the discretely jump-constrained BSDE and the solution 
$(\bar Y^\pi,\bar\Yc^\pi)$ to the above discrete time scheme.

\vspace{1mm}

\begin{Theorem} \label{maintheo}
There exists some constant $C$ such that: 
\beqs
  \E \Big[ \big|Y_{t_k}^\pi - \bar Y_{t_k}^\pi\big|^2 \Big]  
 +   \sup_{t \in (t_k,t_{k+1}]} \E \Big[  \big|Y_t^\pi - \bar Y_{t_{k+1}}^\pi\big|^2 \Big] 
 +    \sup_{t \in [t_k,t_{k+1})} \E \Big[  \big|Y_t^\pi - \bar\Yc_{t_{k}}^\pi\big|^2 \Big]    & \leq & C |\pi|,
\enqs
for all $k$ $=$ $0,\ldots,n-1$. 
\end{Theorem}

\vspace{3mm}

The above  convergence rate $|\pi|^{1\over 2}$ in the $L^2-$norm  for the discretization of the discretely jump-constrained BSDE is the same as for standard BSDE, see \cite{boutou04}, \cite{zha04}. 
By combining with  the convergence result in Section 4,  we finally obtain an estimate on the error due to the discrete time approximation of the minimal solution $Y$ 
to the BSDE with nonpositive jumps. We split the error between the positive and negative parts:
\beqs
{\rm Err}_+^\pi(Y) &:= &   \max_{k\leq n-1} 
\Big(  \E \Big[ \big(Y_{t_k} - \bar Y_{t_k}^\pi\big)_+^2\Big]   +  \sup_{t \in (t_k,t_{k+1}]} \E \Big[ \big(Y_t - \bar Y_{t_{k+1}}^\pi\big)_+^2 \Big] 
+  \sup_{t \in [t_k,t_{k+1})} \E \Big[ \big(Y_t - \bar\Yc_{t_{k}}^\pi\big)_+^2 \Big] \Big)^{1\over 2} \\
{\rm Err}_-^\pi(Y) &:= &   \max_{k\leq n-1} 
\Big(  \E \Big[ \big(Y_{t_k} - \bar Y_{t_k}^\pi\big)_-^2\Big]   +  \sup_{t \in (t_k,t_{k+1}]} \E \Big[ \big(Y_t - \bar Y_{t_{k+1}}^\pi\big)_-^2 \Big] 
+  \sup_{t \in [t_k,t_{k+1})} \E \Big[ \big(Y_t - \bar\Yc_{t_{k}}^\pi\big)_-^2 \Big] \Big)^{1\over 2}.
\enqs

\begin{Corollary}  \label{corollY}
We have:
\beqs
{\rm Err}_-^\pi(Y)   & \leq & C |\pi|^{1\over 2}. 
\enqs
Moreover, under  {\bf (H1')} and {\bf (H2')}, 
\beqs
{\rm Err}_+^\pi(Y)   & \leq & C |\pi|^{1\over 10}, 
\enqs
and when $f(x,a)$ does not depend on $y$:
\beqs
{\rm Err}_+^\pi(Y)   & \leq & C |\pi|^{1\over 6}.  
\enqs
\end{Corollary}
{\bf Proof.} Recall from Proposition \ref{propconv1} that $\Yc_t^\pi$ $\leq$ $Y_t^\pi$ $\leq$ $Y_t$, $0\leq t\leq T$. Then, we have: 
$(Y_{t_k}-\bar Y_{t_k}^\pi)_-$ $\leq$ $|Y_{t_k}^\pi - \bar Y_{t_k}^\pi|$, $(Y_t-\bar Y_{t_{k+1}}^\pi)_-$ $\leq$ $|Y_t^\pi-\bar Y_{t_{k+1}}^\pi|$, and 
$(Y_{t_k}-\bar\Yc_{t_k}^\pi)_-$ $\leq$ $|Y_{t_k}^\pi - \bar\Yc_{t_k}^\pi|$, for all $k$ $\leq$ $n-1$, and $t$ $\in$ $[0,T]$.  The error bound on ${\rm Err}_-^\pi(Y)$ follows then from 
the estimation in Theorem \ref{maintheo}.   The error bound on  ${\rm Err}_-^\pi(Y)$  follows from Corollary \ref{cormainTH} and Theorem \ref{maintheo}. 
\ep

\vspace{5mm}

\begin{Remark}
{\rm  In the particular case where $f$ depends only on $(x,a)$, our discrete time approximation scheme is a probabilistic scheme for the fully nonlinear HJB equation associated to the stochastic control problem \reff{constov}.   As in \cite{kr00}, 
\cite{barjac07} or \cite{ftw11},  we have non symmetric  bounds on the rate of convergence.  For instance, in \cite{ftw11}, the authors  obtained a convergence rate $|\pi|^{1 \over 4}$ on one side and 
$|\pi|^{1\over 10}$ on the other side,  while we improve the rate  to $|\pi|^{1\over 2}$ for one side, and $|\pi|^{1\over 6}$ on  the other side.  This induces a global error  
${\rm Err}^\pi(Y)$ $=$  ${\rm Err}_+^\pi(Y)$ $+$ ${\rm Err}_-^\pi(Y)$ of order $|\pi|^{1\over 6}$, which is derived without any non degeneracy condition on the controlled diffusion coefficient. 
\ep
}
\end{Remark}

\vspace{5mm}

\ni {\bf Proof of Theorem \ref{maintheo}.} 

\vspace{1mm}

\ni Let us  introduce  the continuous time version of \reff{schemeBSDE2}.   
By the martingale representation theorem, there exists $\tilde\Zc^\pi$ $\in$ $L^2(W)$ and $\tilde\Uc^\pi$ $\in$ $L^2(\tilde\mu)$ such that
\beqs
\bar Y_{t_{k+1}}^\pi &=& \E_{t_k} \big[ \bar Y_{t_{k+1}}^\pi \big] + \int_{t_k}^{t_{k+1}} \tilde\Zc^\pi_t dW_t + \int_{t_k}^{t_{k+1}} \int_A \tilde\Uc_t^\pi(a) \tilde\mu(dt,da), \;\; k < n, 
\enqs
and we can then define the continuous-time  processes $\bar Y^\pi$ and $\bar\Yc^\pi$ by: 
\beq
\bar \Yc ^\pi_{t} & = & \bar Y ^\pi_{t_{k+1}} +(t_{k+1}-t)f(\bar X ^\pi _{t_k}, I_{t_k},\bar\Yc_{t_k}^\pi ) \label{Ycpicont}   \\
& & \;\;\; - \int_t^{t_{k+1}} \tilde \Zc ^\pi_{t}dW_t - \int_{t}^{t_{k+1}} \int_A \tilde\Uc_t^\pi(a) \tilde\mu(dt,da),  \;\;\;\;\;\;\;  t\in[t_k,t_{k+1}),\nonumber \\
\bar Y ^\pi_{t} & = & \bar Y ^\pi_{t_{k+1}} +(t_{k+1}-t)f(\bar X ^\pi _{t_k}, I_{t_k},\bar\Yc_{t_k}^\pi ) \label{Ypicont}  \\
& & \;\;\; -\int_t^{t_{k+1}} \tilde \Zc ^\pi_{t}dW_t - \int_{t}^{t_{k+1}} \int_A \tilde \Uc_t^\pi(a) \tilde\mu(dt,da),  \;\;\;\;\;\;\;  t\in (t_k,t_{k+1}], \nonumber 
\enq
for $k=0,\ldots,n-1$. 
Denote by $\delta Y_t^\pi$ $=$ $Y_t^\pi-\bar Y_t^\pi$, $\delta \Yc_t^\pi $ $=$ $\Yc_t^\pi-\bar \Yc_t^\pi$,  $\delta \Zc_t^\pi$ $=$ $\Zc_t^\pi-\tilde \Zc_t^\pi$,  $\delta \Uc_t^\pi$ 
$=$ $\Uc_t^\pi-\tilde \Uc_t^\pi$ and  
$\delta f_t$ $=$ $f(X_t,I_t,\Yc_t^\pi)-f(\bar X^\pi_{t_k}, I_{t_k},\bar\Yc_{t_k}^\pi )$ for $t\in[t_k,t_{k+1})$. Recalling \reff{BSDEDC2} and \reff{Ycpicont}, we have by It\^o's formula: 
\beqs
\Delta_ t & := & \E_{t_k} \Big[ |\delta \Yc_t^\pi|^2  +  \int_t^{t_{k+1}} |\delta \Zc_s^\pi|^2 ds  + \int_t^{t_{k+1}} \int_A |\delta \Uc_s^\pi(a)|^2 \lambda(da) ds \Big] \\
&=&  \E_{t_k} \big[ |\delta Y_{t_{k+1}}^\pi|^2\big| \big] +  \E_{t_k}\Big[ \int_t^{t_{k+1}} 2\delta\Yc_s^\pi\delta f_s \Big] ds 
\enqs
for all $t\in[t_k,t_{k+1})$.  By the Lipschitz continuity of $f$ in {\bf (H2)} and  Young inequality, we then have:
\beqs
\Delta_t & \leq &  \E_{t_k} \big[ |\delta Y_{t_{k+1}}^\pi|^2\big| \big]  
+ \E_{t_k}\Big[ \int_t^{t_{k+1}}  \eta  |\delta\Yc_s^\pi|^2  ds  + \frac{C}{\eta}  \pi |\delta\Yc_{t_k}^\pi|^2 \Big] \\
& & \; + \;  \frac{C}{\eta}  \E_{t_k}\Big[ \int_t^{t_{k+1}} \big(|X_s-\bar X_{t_k}^\pi|^2 + |I_s -I_{t_k}|^2  + |\Yc_{s}^\pi - \Yc_{t_k}^\pi|^2 \big) ds\Big]. 
\enqs
From Gronwall's  lemma, and by taking $\eta$ large enough,   this yields for all $k$ $\leq$ $n-1$: 
\beq
\E_{t_k} \Big[ |\delta \Yc_{t_k}^\pi|^2  \Big] & \leq & 
e^{C |\pi|} \E_{t_k} \big[ |\delta Y_{t_{k+1}}^\pi|^2\big| \big]  + C B_k  \label{estimYgron1} 
\enq
where
\beq
B_k & =&  \E_{t_k}\Big[ \int_{t_k}^{t_{k+1}} \big( |X_s-\bar X_{t_k}^\pi|^2 + |I_s -I_{t_k}|^2 +  |\Yc_{s}^\pi - \Yc_{t_k}^\pi|^2  \big) ds  \Big]  \nonumber \\
& \leq &  C |\pi| \Big(  \E_{t_k}\big[ \sup_{s \in [t_k,t_{k+1}]} |X_s-\bar X_{t_k}^\pi|^2 \big]  + |\pi|(1 + |X_{t_k}|) \Big), \label{Bestim}
\enq
by \reff{estimI} and Proposition \ref{regulYdis}.  Now, by definition of $Y_{t_{k+1}}^\pi$ and $\bar Y_{t_{k+1}}^\pi$, we have
\beq \label{Ycinter}
|\delta Y_{t_{k+1}}^\pi|^2 & \leq & \esssup_{a\in A} \E_{t_{k+1},a}\big[ | \delta \Yc_{t_{k+1}}^\pi|^2 \big].
\enq
By plugging \reff{Bestim}, \reff{Ycinter} into \reff{estimYgron1},  taking conditional expectation with respect to $I_{t_k}$ $=$ $a$, and taking essential supremum over $a$, we obtain: 
\beqs 
\esssup_{a\in A} \E_{t_k,a} \Big[ |\delta \Yc_{t_k}^\pi|^2  \Big] & \leq & e^{C |\pi|} \esssup_{a\in A} \E_{t_k,a} \big[  \esssup_{a\in A} \E_{t_{k+1},a}\big[ | \delta \Yc_{t_{k+1}}^\pi|^2 \big] \nonumber \\
& & \; + \; C|\pi| \Big( \esssup_{a\in A} \E_{t_k,a}\big[ \sup_{s \in [t_k,t_{k+1}]} |X_s-\bar X_{t_k}^\pi|^2 \big]  + |\pi|(1 + |X_{t_k}|) \Big).  \label{interestim}
\enqs
By taking conditional expectation with respect to $\Fc_{t_{k-1}}$, and $I_{t_{k-1}}$ $=$ $a$, taking essential supremum over $a$ in the above inequality, and iterating this backward procedure until time $t_0$ $=$ $0$,  we obtain: 
\beq
\Ec_k^\pi(\Yc) & \leq & e^{C |\pi|}  \Ec_{k+1}^\pi(\Yc) + C |\pi| \big( \Ec_k^\pi(X) + |\pi|(1 + \E[|X_{t_k}|]) \big) \nonumber \\
& \leq & e^{C |\pi|}  \Ec_{k+1}^\pi(\Yc) + C |\pi|^2,   \;\;\;\;\;   k\leq n-1, \label{estimYkk+1}
\enq
where we recall the auxiliary error control $\Ec_k^\pi(X)$ on $X$ in \reff{errorauxX} and its estimate in Lemma \ref{lemEulerX}, and set: 
\beqs
\Ec_k^\pi(\Yc) &:=& \E\Big[\esssup_{a\in A}\E_{t_1,a}\big[ \ldots \esssup_{a\in A} \E_{t_k,a}\big[ |\delta \Yc_{t_k}^\pi|^2 \big] \ldots \big] \Big]. 
\enqs
By a direct induction on \reff{estimYkk+1}, and recalling that $n|\pi|$ is bounded, we get
\beqs
\Ec_k^\pi(\Yc) & \leq &  C \big(  \Ec_n^\pi(\Yc)  + |\pi| \big)  \\
& \leq & C (\Ec_n^\pi(X) + |\pi| \big) \; \leq \; C |\pi|,
\enqs
since $g$ is Lipschitz, and using again the estimate in Lemma \ref{lemEulerX}. Observing that $\E[|\delta Y_{t_k}^\pi|^2]$,  
$\E[|\delta\Yc_{t_k}^\pi|^2]$ $\leq$ $\Ec_k^\pi(\Yc)$,  we get the estimate: 
\beqs
\max_{k\leq n} \E\big[|Y_{t_k}^\pi-\bar Y_{t_k}^\pi|^2\big]  + \E\big[|\Yc_{t_k}^\pi-\bar\Yc_{t_k}^\pi|^2\big] & \leq & C |\pi|. 
\enqs
Moreover, by Proposition \ref{regulYdis}, we have
\beqs
 \sup_{t\in[t_k,t_{k+1})} \E \Big[ |\Yc_t^\pi - \Yc_{t_{k}}^\pi|^2 \Big]  
+ \sup_{t\in(t_k,t_{k+1}]}  \E \Big[ |Y_t^\pi - Y_{t_{k+1}}^\pi|^2 \Big]  & \leq & C(1+ \E[|X_{t_k}|]) |\pi| \\
& \leq & C(1+|X_0|)|\pi|. 
\enqs
This implies finally that:
\beqs
 \sup_{s\in(t_k,t_{k+1}]} \E\Big[  |Y_t^\pi - \bar Y_{t_{k+1}}^\pi|^2 \Big] & \leq &   
 2 \sup_{s\in(t_k,t_{k+1}]}  \E \Big[ |Y_t^\pi - Y_{t_{k+1}}^\pi|^2 \Big] + 2  \E \Big[ |Y_{t_{k+1}}^\pi - \bar Y_{t_{k+1}}^\pi|^2 \Big] \\
 & \leq & C|\pi|, 
\enqs
as well as
\beqs
 \sup_{s\in[t_k,t_{k+1})} \E\Big[  |Y_t^\pi - \bar\Yc_{t_{k}}^\pi|^2 \Big] & \leq &   
 2 \sup_{s\in[t_k,t_{k+1})}  \E \Big[ |Y_t^\pi - \Yc_{t_{k}}^\pi|^2 \Big] + 2  \E \Big[ |\Yc_{t_{k}}^\pi - \bar\Yc_{t_{k}}^\pi|^2 \Big] \\
 & \leq & C|\pi|. 
\enqs
\ep

\subsection{Approximate optimal control} \label{approxcontrol}

 We now consider the special  case  where $f(x,a)$ does not depend on $y$, so that  the discrete time scheme \reff{schemeBSDE} is an approximation for the value function of the stochastic control problem:
\beq 
V_0  & :=& \sup_{\alpha \in \Ac}  J(\alpha)  \; = \; Y_0, \label{stoV0} \\
J(\alpha) & = & \E \Big[ \int_0^T f(X_t^{\alpha},\alpha_t) dt + g(X_T^{\alpha}) \Big], \nonumber 
\enq
where $\Ac$ is the set of $\G$-adapted  control processes $\alpha$ valued in $A$, and $X^\alpha$ is the  controlled diffusion in $\R^d$: 
\beqs
X_t^\alpha  &=& X_0 + \int_0^t b(X_s^\alpha,\alpha_s) ds +  \int_0^t \sigma(X_s^\alpha,\alpha_s) dW_s, \;\;\; 0 \leq t\leq T. 
\enqs 
(Here $\G$ $=$ $(\Gc_t)_{0\leq t\leq T}$ denotes some filtration under which $W$ is a standard Brownian motion).  
Let us now define the discrete time version of  \reff{stoV0}. We introduce the set $\Ac^\pi$ of discrete time  processes $\alpha$ $=$ $(\alpha_{t_k})_k$ with $\alpha_{t_k}$ $\Gc_{t_k}$-measurable, and valued in 
$A$. For each $\alpha$ $\in$ $\Ac^\pi$, we consider the controlled discrete time process $(X_{t_k}^{\pi,\alpha})_k$ of Euler type defined by:
\beqs
X_{t_k}^{\pi,\alpha} &=&  X_0 + \sum_{j=0}^{k-1}  b(X_{t_j}^{\pi,\alpha},\alpha_{t_j}) \Delta t_j +   \sum_{j=0}^{k-1} \sigma(X_{t_j}^{\pi,\alpha},\alpha_{t_j}) \Delta W_{t_j},  \;\;\;  k  \leq n,
\enqs
where $\Delta W_{t_j}$ $=$ $W_{t_{j+1}}-W_{t_j}$, and the gain functional:
\beqs
J^\pi(\alpha) &=& \E \Big[ \sum_{k=0}^{n-1} f(X_{t_k}^{\pi,\alpha},\alpha_{t_k}) \Delta t_k + g(X_{t_n}^{\pi,\alpha}) \Big].  
\enqs
Given any $\alpha$ $\in$ $\Ac^\pi$, we  define its continuous time piecewise-constant interpolation $\alpha$ $\in$ $\Ac$  by setting: $\alpha_t$ $=$ $\alpha_{t_k}$, for $t\in[t_k,t_{k+1})$ (by misuse of notation, we keep the same notation $\alpha$ for the discrete time and continuous time  interpolation). By standard arguments similar  to those  for Euler scheme of SDE,  there exists some positive constant $C$ such that for all 
$\alpha$ $\in$ $\Ac^\pi$,  $k$ $\leq$ $n-1$:
\beqs
\E\Big[ \sup_{t \in [t_k,t_{k+1}]} \big| X_t^\alpha - X_{t_k}^{\pi,\alpha} \big|^2 \Big] & \leq & C |\pi|, 
\enqs
 from which we easily deduce by Lipschitz property of $f$ and $g$: 
 \beq \label{Jpiestim}
 \big| J(\alpha) - J^\pi(\alpha) | & \leq & C |\pi|^{1\over 2}, \;\;\; \forall \alpha \in \Ac^\pi. 
 \enq
 
Let us now consider at each time step $k$ $\leq$ $n-1$, the function $\hat a_k(x)$ which attains the supremum over $a$ $\in$ $A$ of $\bar\vartheta_k^\pi(x,a)$ in the scheme \reff{schemeBSDE3},  so that: 
\beqs
\bar v_k^\pi(x) & = & \bar\vartheta_k^\pi\big(x,\hat a_k(x)\big), \;\;\; k=0,\ldots,n-1. 
\enqs 
Let us define the process $(\hat X_{t_k}^\pi)_k$ by: $\hat X_0^\pi$ $=$ $X_0$, 
\beqs
\hat X_{t_{k+1}}^\pi &=& \hat X_{t_k}^\pi + b(\hat X_{t_k}^\pi,\hat a_k(\hat X_{t_k}^\pi)) \Delta t_k + \sigma(\hat X_{t_k}^\pi,\hat a_k(\hat X_{t_k}^\pi)) \Delta W_{t_k}, \;\;\; k \leq n-1, 
\enqs
and notice  that  $\hat X^\pi$ $=$ $X^{\pi,\hat\alpha}$, where $\hat\alpha$ $\in$ $\Ac^\pi$ is a feedback control defined  by:
\beqs
\hat\alpha_{t_k} &=& \hat a_k(\hat X_{t_k}^\pi) \; = \; \hat a_k(X_{t_k}^{\pi,\hat\alpha}) , \;\;\; k =0,\ldots,n. 
\enqs
Next, we observe that the conditional  law of $\bar X_{t_{k+1}}^\pi$ given $(\bar X_{t_k}^\pi=x, I_{t_k}=\hat a_k(\bar X_{t_k}^\pi)=\hat a_k(x))$ is the same than the conditional law of $X_{t_{k+1}}^{\pi,\hat\alpha}$  given 
$X_{t_k}^{\pi,\hat\alpha}$ $=$  $x$,  for $k$ $\leq$ $n-1$, and thus the induction step in the scheme \reff{schemeBSDE2} or \reff{schemeBSDE3}  reads as: 
\beqs
\bar v_k^\pi(X_{t_k}^{\pi,\hat\alpha}) &=& \E \Big[ \bar v_{k+1}^\pi(X_{t_{k+1}}^{\pi,\hat\alpha}) \big| X_{t_k}^{\pi,\hat\alpha} \Big] +  f(X_{t_k}^{\pi,\hat\alpha},\hat\alpha_{t_k}) \Delta t_k,  
\;\;\; k\leq n-1. 
\enqs
By induction, and law of  iterated conditional expectations, we then get:
\beq \label{Y0alpha}
\bar Y_0^\pi \; = \;  \bar v_0^\pi(X_0) &=& J^\pi(\hat\alpha). 
\enq
Consider the continuous time piecewise-constant interpolation $\hat\alpha$ $\in$ $\Ac$  defined by: $\hat\alpha_t$ $=$ $\hat\alpha_{t_k}$, for $t\in[t_k,t_{k+1})$.  By \reff{stoV0}, 
\reff{Jpiestim}, \reff{Y0alpha}, and Corollary \ref{corollY},  we finally obtain:
\beqs
0   \; \leq \; V_0 - J(\hat\alpha) & = & Y_0 - \bar Y_0^\pi +  J^\pi(\hat\alpha) - J(\hat\alpha) \\
& \leq & C|\pi|^{1\over 6} + C |\pi|^{1 \over 2} \; \leq \; C|\pi|^{1\over 6},
\enqs
for $|\pi|$ $\leq$ $1$.  In other words,  for any small $\eps$ $>$ $0$, we  obtain  an $\eps$-approximate optimal control $\hat\alpha$ for the stochastic control problem \reff{stoV0} by taking $|\pi|$ of order  $\eps^{6}$.

 \vspace{3mm}


\end{document}